\documentclass[12pt,twoside]{article}
\usepackage{epsfig}
\usepackage{latexsym}
\usepackage{graphicx}
\usepackage{arabtex}
 \def\Dj{\hbox{D\kern-.73em\raise.30ex\hbox{-} \raise-.30ex\hbox{}}}
 \def\dj{\hbox{d\kern-.33em\raise.80ex\hbox{-} \raise-.80ex\hbox{\kern-.40em}}}

\vspace{10mm}

\begin{document}

\vspace{10mm}

\baselineskip=0.20in
\begin{center}
\textbf{A SURVEY ON THE SPECIAL FUNCTION }
{\includegraphics[width=0.36in,height=0.39in]{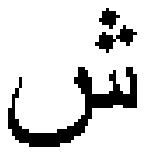}}
\end{center}

\vspace{10mm}

\begin{center}
\textbf{ ANDREA OSSICINI }
\end{center}

\vspace{20mm}

\begin{center}
\begin{figure}[htbp]
\centerline{\includegraphics[width=1.44in,height=1.44in]{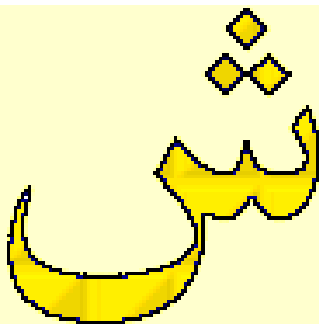}}
\label{fig1}
\end{figure}

\end{center}
\setarab

\baselineskip=0.30in

\setnashbf

\vspace{20mm}

\noindent {\small {\bf Abstract.} The purpose of the work is to
furnish a complete study of a discrete and special function,
discovered by the author and named with the Arabian letter <^s>
(Shin)\footnote{ The letter <^s> is the thirteenth letter of the
Arabian alphabet.}.

It includes three other papers, published in the international
journal "\textbf{Kragujevac Journal of Mathematics}".

The methods, the techniques and the style of the demonstration
inside such a work, are all related to \textbf{Leonhard Euler}, the
very great Swiss mathematician.

\newpage

 \setcounter{page}{01}
 \baselineskip=0.20in
\textit{Kragujevac J. Math.} 27 (2005) 01--25.

\vspace{20mm}

\setarab

\baselineskip=0.30in

\setnashbf
\begin{center}
\textbf{THE SPECIAL FUNCTION } <^s> \textbf{.}
\end{center}
\setnash
\begin{center}
\vspace{10mm}

{\large \bf  Andrea Ossicini}

\vspace{10mm}

\baselineskip=0.20in

{\it Via delle Azzorre 352-D2, 00121 Roma, Italy\\
{\rm (e-mail: a.ossicini@finsiel.it)}

\vspace{10mm}

(Received March 12, 2005)
}
\end{center}
\vspace{10mm}

\baselineskip=0.17in \noindent {\small {\bf Abstract.} The purpose
of the work is to furnish a study of a discrete and special
function, discovered by the author and named with the Arabian letter
<^s> (Shin)\footnote{ The letter <^s> is the thirteenth letter of
the Arabian alphabet.}, produced by a family of functions. A
fundamental theorem is enunciated: it connects the rational values
of this family with a natural number; to this aim two rational-value
functions will be created, with the characteristic, in the field of
real positive numbers, to be piecewise continuous. We state not only
the existence of a separation element, but we prove that this
element is just formed by one only integer constant function, the
value of which is equal to 2. We point out a hypothetical, subtle
connection among the special function <^s> , the Eulerian function
Gamma and the second-order Eulerian numbers. It is finally proved
that <^s> is completely monotonic: this characteristic is peculiar
for the functions that have considerable applications in different
fields of pure and applied Mathematics.

\vspace{6mm}

2000\textit{ Mathematics Subject Classification: primary
}33B15\textit{; secondary }26A48.

\vspace{5mm} Keywords: special functions, gamma function,
hypergeometric series, second-order Eulerian numbers, completely
monotonic functions, integral transforms.

\baselineskip=0.20in
\newpage
\begin{center}
1. FROM THE CONTINUOUS TO THE DISCRETE
\end{center}

\vspace{4mm}

Let's consider the three \textit{following} transcendental
functions, determined by the letter <^s> : \linebreak

\qquad <^s> $\left[1\right]=\left({{\kern
2pt}\,1+\frac{1}{3k}\;}\right)^{^{2k+1}}\;;$ <^s>
$\left[2\right]=\left({{\kern
2pt}\,1+\frac{1}{3k-1}\;}\right)^{^{2k+1}}\;$; <^s>
$\left[3\right]=\left({{\kern
2pt}\,1+\frac{1}{3k-2}\;}\right)^{^{2k+1}}\;\;$

\vspace{4mm}

\noindent and where the variable $k$ will be restricted only in
the real positive values, except some values, explained later.

The analytic study of the three functions permits to identify their
progress, for the presence, in each of them, of two asymptotes: a horizontal
asymptote, got by the limit for $k$ going to infinity, represented by the
horizontal straight line, of height\footnote{ $e$ represents the
Euler's number.}:
\[e^{2 \mathord{\left/ {\vphantom {2 3}} \right.
\kern-\nulldelimiterspace} 3} \approx \;1,947734041\]

A vertical asymptote, specific for every function, got by the
following limits:

\vspace{4mm}

\[ \lim_{k\to 0^+}  \ \mbox{<^s>}  \
\left[ 1 \right]\;=\infty \mbox{ ; }\lim_{k\to \frac{1}{3}^+}  \
\mbox{<^s>}  \ \left[ 2 \right]\;=\infty \mbox{ ; }\lim_{k\to
\frac{2}{3}^+}  \ \mbox{<^s>} \ \left[ 3 \right]\;=\infty \]

\vspace{4mm}

Besides by the calculus of the first derivative of each function and
the study of its sign, it's possible to verify that the <^s>
functions are decreasing in their whole field of existence, more
precisely, if we consider for the variable $k$ the whole positive
real axis, we must exclude at least for the second and third
function, respectively the intervals $\left(
{0,\raise0.7ex\hbox{$1$} \!\mathord{\left/ {\vphantom {1
3}}\right.\kern-\nulldelimiterspace}\!\lower0.7ex\hbox{$3$}}
\right]$ and $\left( {0,\raise0.7ex\hbox{$2$} \!\mathord{\left/
{\vphantom {2
3}}\right.\kern-\nulldelimiterspace}\!\lower0.7ex\hbox{$3$}}
\right]$. \vspace{4mm}

In fact, if we derive the <^s> functions we get respectively:

\vspace{3mm}

$\frac{d}{dk} \mbox{<^s>} \left[ 1 \right]\;=\mbox{<^s>} \left[ 1
\right]\;\;\cdot \left[ {2\cdot \log
\;(1+\frac{1}{3k})-\frac{2k+1}{k(3k+1)}} \right]\prec 0$ $\forall
\;\;k\in \mathbf{R^{+}}$

\vspace{3mm}

$\frac{d}{dk} \mbox{<^s>} \left[ 2 \right]\;=\mbox{<^s>} \left[ 2
\right]\;\;\;\cdot \left[ {2\cdot \log
\;(1+\frac{1}{3k-1})-\frac{2k+1}{k(3k-1)}} \right]\prec 0 \quad
\forall \;\;k\in \mathbf{R^{+}}\quad$ and
$\;\;k>\raise0.7ex\hbox{$1$} \!\mathord{\left/ {\vphantom {1
3}}\right.\kern-\nulldelimiterspace}\!\lower0.7ex\hbox{$3$}$

\vspace{3mm}

$\frac{d}{dk} \mbox{<^s>} \left[ 3 \right]\;=\mbox{<^s>} \left[ 3
\right]\;\cdot \;\left[ {2\cdot \log
\;(1+\frac{1}{3k-2})-\frac{3\cdot \left( {2k+1} \right)}{\left(
{3k-1} \right)\cdot \left( {3k-2} \right)}} \right]\prec 0$ $\forall
\;\;k\in \mathbf{R^{+}}\quad$ and $\;\;k>\raise0.7ex\hbox{$2$}
\!\mathord{\left/ {\vphantom {2
3}}\right.\kern-\nulldelimiterspace}\!\lower0.7ex\hbox{$3$}$

\vspace{3mm}

In the Fig. 1 we have represented three functions where, among other things,
it's evident that all of them have only one point of intersection with the
horizontal straight line of height 2.

\medskip
\centerline{\epsfig{figure=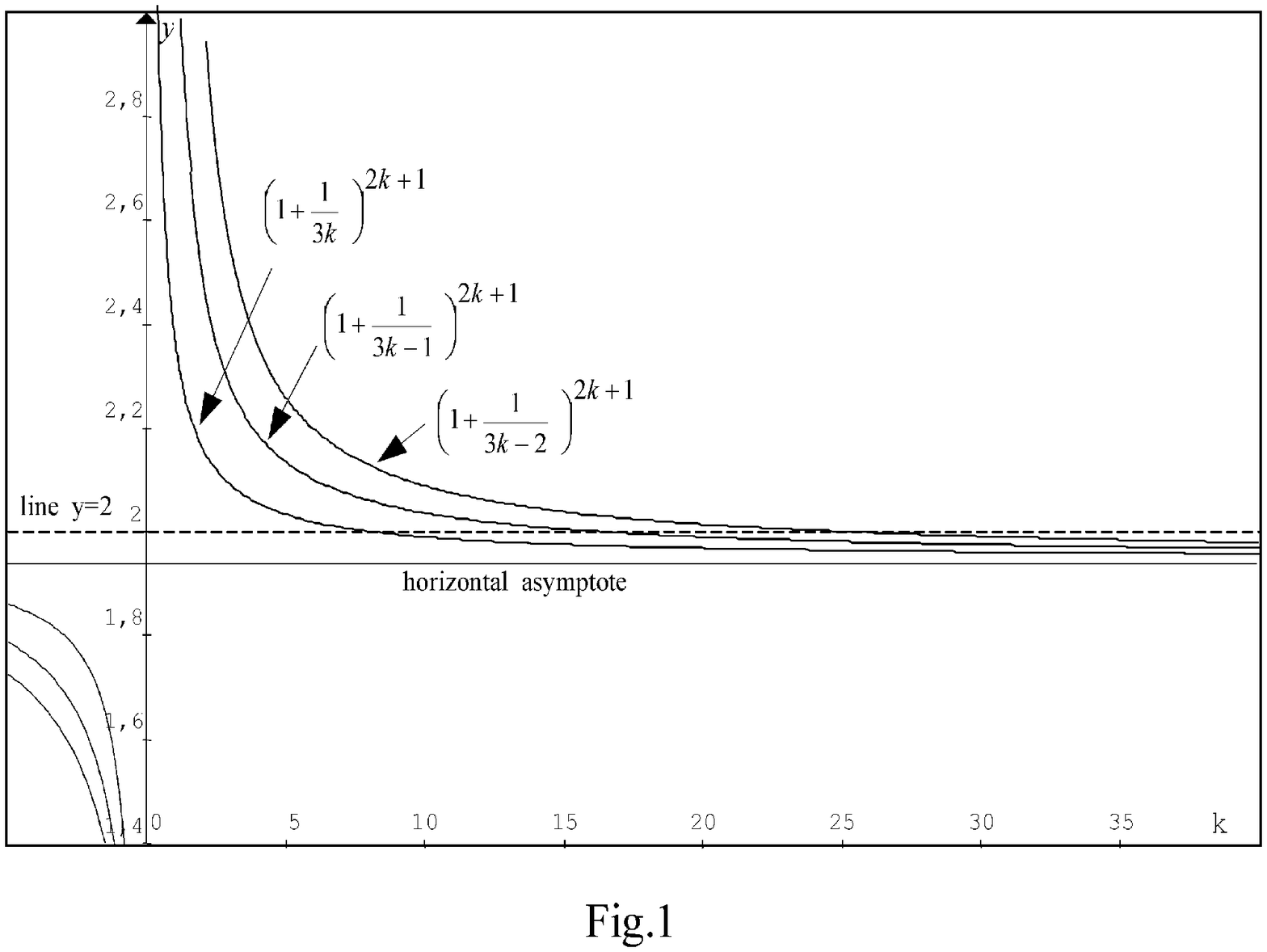,width=14cm,height=8cm}}
\medskip
\baselineskip=0.20in

That being stated, let's proceed in the passage from the
continuous to the discrete, by considering for $k$, the only
integer positive values; under these hypotheses it's possible to
verify:

\vspace{3mm}

\qquad \qquad $\mbox{<^s>} \left[ 1 \right]\;  \;\succ \;2$ for
$1\le k\le 8$ \qquad and \qquad $\mbox{<^s>} \left[ 1 \right]\;
\prec \;2$ for $k\ge 9$

\qquad \qquad $\mbox{<^s>} \left[ 2 \right]\;  \;\succ \;2$ for
$9\le k\le 16$ \qquad and \qquad $\mbox{<^s>} \left[ 2 \right]\;
\prec \;2$ for $k\ge 17$

\qquad \qquad $\mbox{<^s>} \left[ 3 \right]\; \succ \;2$ for $17\le
k\le 25$ \qquad and \qquad $\mbox{<^s>} \left[ 3 \right]\; \prec 2$
for $k\ge 26$

\vspace{3mm}

\noindent and therefore, in the discrete it's possible to define
some intervals\footnote{ $\ell \ $shows the ``ordinal number'' of
the interval $I_\ell $.} $I_\ell$ of integer values of the
variable $k$, in order to characterize some limitations of the
values that have the above stated functions; in fact it's possible
to verify that for appropriate intervals $I_\ell ;\ell =2,3$ the
following are valid (see Fig. 2):

\vspace{3mm}

$\mbox{<^s>} \left[ 1 \right]  \prec 2\prec   \mbox{<^s>} \left[ 2
\right]\;$ for $9\le k\le 16 \Rightarrow I_2 $; $\mbox{<^s>} \left[
2 \right] \prec 2\prec   \mbox{<^s>} \left[ 3 \right]\;$ for $17\le
k\le 25\Rightarrow I_{3} $

\medskip
\centerline{\epsfig{figure=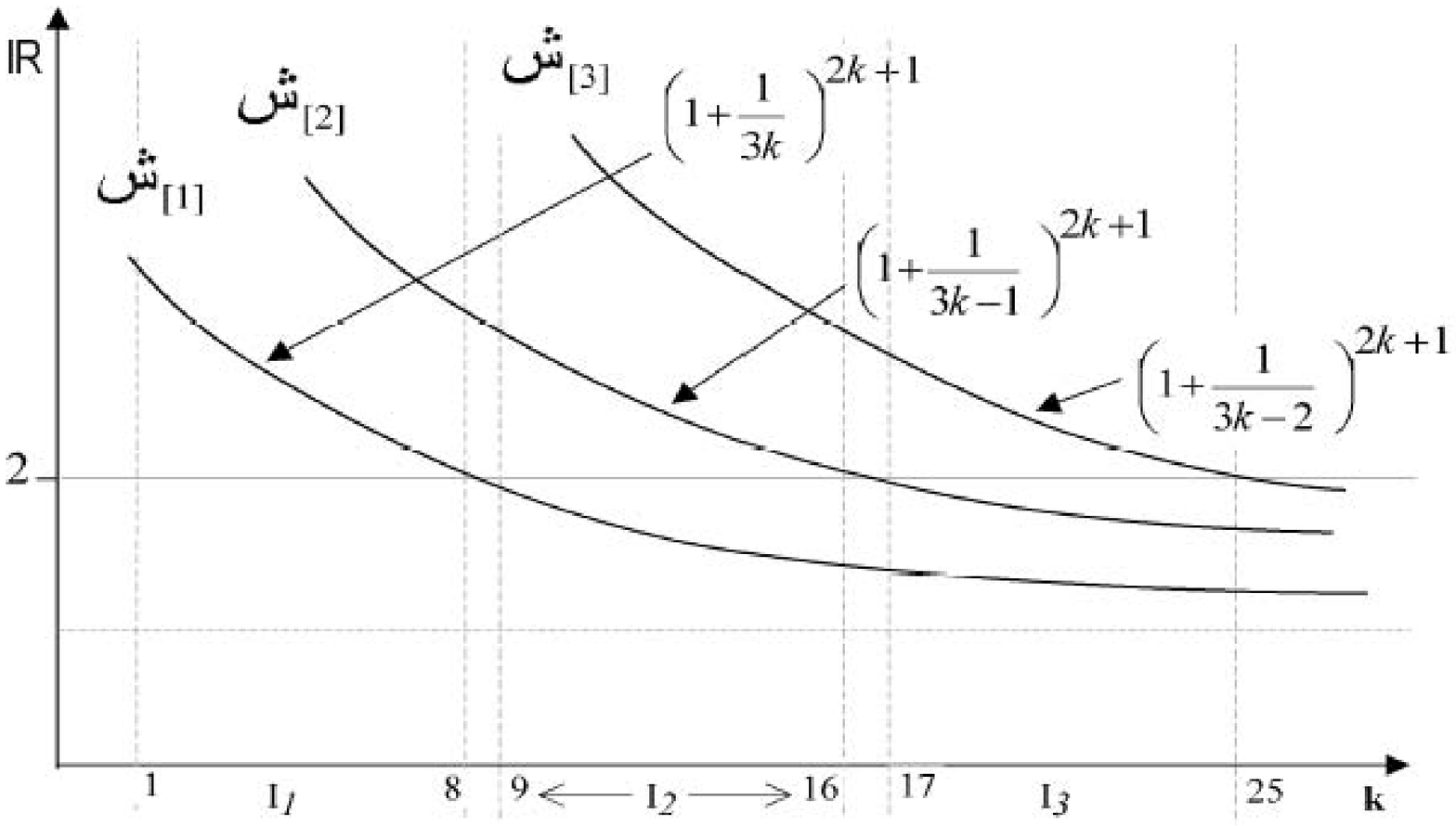,width=16cm,height=9cm}}
\medskip
\baselineskip=0.20in
\begin{center}
{\bf Fig. 2}
\end{center}

\vspace{4mm}

In conclusion it's possible to introduce a family of <^s> functions,
by the definition of appropriate arcs, whose separation element is
the horizontal straight line of height 2.

The construction-algorithm of the above stated family is therefore
describable in the following way: we begin from the first algebraic
expression of the <^s> function, that is $\;\left(
{\;1+\frac{1}{3k}\;} \right)^{\;2k+1}\;$, starting to calculate by
growing values of the integer positive variable $k$, the
corresponding rational values of $\mbox{<^s>} \left[ 1 \right]\;$ ;
for the first 8 integer values of $k$, the function has rational
values greater than 2 and it's therefore possible to associate to
such values a bounded arc of the same function, represented in the
discrete field by a sequence of rational numbers, each of them
greater than 2.

After that we decrease of a unity the value of the denominator of
the fraction inside the <^s> function, consequently we'll get the
algebraic expression of a new function, that is $\;\left(
{1+\frac{1}{3k-1}} \right)^{2k+1}\;$, that we have previously
identified with $\mbox{<^s>} \left[ 2 \right]\;$ and that we can
define as the \textit{following} of $\mbox{<^s>} \left[ 1
\right]\;$.
\newpage
We'll repeat, for it too, the same procedure and therefore we'll
calculate by growing values of the integer variable $k$, but greater
than the previous 8,
 the corresponding rational values of $\mbox{<^s>} \left[ 2 \right]\;$;
also in this case for 8 integer values of $k$ , the function has
rational values greater than 2 and therefore it's possible to
associate to such values a bounded arc of the same function,
represented in the discrete field by a sequence of rational numbers,
everyone greater than 2.

Now, for the same 8 integer values of $k$, utilized for $\mbox{<^s>}
\left[ 2 \right]\;$, it's besides possible to verify that the
$\mbox{<^s>} \left[ 1 \right]\;$ function has, on the contrary,
rational values smaller than 2 and it's therefore possible to
associate to such values a bounded arc of the same function,
represented in the discrete field by a sequence of rational numbers,
all smaller than 2.

If we repeat the procedure and then we decrease, as usual, the value
of the denominator of the fraction inside the <^s> function, we can
build the $\mbox{<^s>} \left[ 3 \right]\;$ function, that so results
the \textit{following} of  $\mbox{<^s>} \left[ 2\right]\;$.

In this case, differently from the first two, exactly for 9 integer values
of $k$, greater than the previous 8, the function has rational values greater
than 2 and therefore it's possible to associate to such values a bounded arc
of the same function, represented in the discrete field by a sequence of
rational numbers, all of them greater than 2.

Here, for the same 9 integer values of $k$, it's possible to verify
that the $\mbox{<^s>}\left[ 2 \right]\;$ function has, on the
contrary, rational values smaller than 2, too, and it's so possible
to associate to such values a bounded arc of the same function,
represented, in the discrete field, by a sequence of rational
numbers, each of them smaller than 2.

Consequently about what described, if we consider the second
interval of 8 integer values of $k$, that is
$k=$9,10,11,12,13,14,15,16, we can build two arcs, represented by
two sequences of rational numbers: the first sequence formed by
numbers greater than 2, because belonging to the arc of the
$\mbox{<^s>} \left[2 \right]\;$ function, the second sequence
formed by numbers smaller than 2, because belonging to the arc of
the $\mbox{<^s>} \left[ 1 \right]\;$ function (see Fig. 2).

What shown is repeatable and it's possible to experiment, while the
integer variable $k$ grows, the determination of two appropriate
arcs, belonging to two \textit{following} <^s> functions.

To simplify the use and the control of the described algorithm a
\textit{vector} function has been defined by the software product
\textit{DERIVE}\footnote{ \textit{DERIVE} is a powerful instrument
of CAS (Computer Algebra System ), spread by Texas
Instruments.}\textit{ Version 6 } for WINDOWS, which allows the
display of 11 consecutive values of a generic <^s> function, by two
only parameters: an integer value of $k$ and a further integer
value, corresponding to the value of the interval (decreased of a
unity) that we wish to study.

In APPENDIX it's given, besides the macro function, which identifies
the vector function, the result of a display got by its use.
\newpage

\setnashbf
\begin{center}
2. THE FUNDAMENTAL THEOREM OF THE <^s> FUNCTION
\end{center}

\setnash

Generalizing what described in the first paragraph and using the
same method to build a couple of <^s> functions, identifiable by a
precise interval $I_\ell $, it's possible therefore to enunciate the
following \textbf{Fundamental Theorem: }

\vspace{1mm}

\textit{Let's consider }$k,\ell $\textit{ natural numbers different from
  zero,} $I_\ell $\textit{ an interval of integer values of the variable k,
  }$\Omega(I_\ell )$\textit{an auxiliary integer function, then for
  each interval }$I_\ell $\textit{ it's always possible the construction of
  a couple of functions with rational values, exclusively depending on k
  and that we denote with }<^s>\textit{, so as the following boundary is
always valid}\footnote{ The boundary can include the sign ``='' if the
  integer variable $k$ goes towards the infinity.}:

\vspace{3mm}

$\mbox{<^s>} \left[ {k,\Omega \;(I_\ell -1)} \right] \prec 2\prec
 \mbox{<^s>} \left[ {k,\Omega \;(I_\ell )} \right]\;\;$ with
$\Omega \;(I_\ell -1)=\Omega \;(I_\ell )-1$ ; $\forall \;I_\ell
\;,\;k,\ell \in \mathbf{N}$

\vspace{3mm}

\noindent \textit{and for the generic }<^s>\textit{ function it's
valid:}

\vspace{1mm}

\begin{equation}
\label{eq1} \mbox{<^s>} \left[ {k,\Omega \;(I_\ell )} \right]\;\quad
= \;\left( {1+\frac{1}{3k-\Omega \;(I_\ell )}} \right)^{2k+1}\;
\end{equation}

The \textit{auxiliary }integer function $\Omega(I_\ell )$, that really
represents a growing ``step function'', is defined, for the intervals of 8
or 9 following values of $k$, in the following way:

\begin{itemize}
\item $\Omega \;(I_1 )=0$ for $k$=1,{\ldots},8
\item $\Omega \;(I_2 )=1$ for $k$=9,{\ldots},16 ; $\Omega \;(I_3 )$=2 for $k$=17,{\ldots},25 ; $\Omega \;(I_4 )$=3 for $k$=26,{\ldots},34
\item $\Omega \;(I_5 )=4$ for $k$=35,{\ldots},43 ; $\Omega \;(I_6 )$ =5 for $k$=44,{\ldots},51 ; $\Omega \;(I_7 )$=6 for $k$=52,{\ldots},60
\item $\Omega \;(I_8 )=7$ for $k$=61,{\ldots},69; $\Omega \;(I_9 )$ =8 for $k$=70,{\ldots},78 ; $\Omega \;(I_{10} )$=9 for $k$=79,{\ldots},86
\item $\Omega \;(I_{11} )=10$ for $k$=87,{\ldots},95 ; $\Omega \;(I_{12} )$=11 for $k$=96,{\ldots},104. ; etc.
\end{itemize}

The sequence of growing values, defined for the \textit{auxiliary}
\footnote{ The function $\Omega(I_\ell )$ is a depending function on
  the integer variable $k$.}function $\Omega(I_\ell )$, is therefore
obtained to allow the individuation of two bounded arcs, belonging
to two \textit{following} <^s> functions, one above the line of
height 2 and the other below it.

The set of the 11 specific intervals, by which the function
$\Omega(I_\ell )$ is always positive and growing, represents a
dominant characteristic of the family of <^s>  functions.

For a better precision we write the \textit{numerical series} that
identifies the extent of such intervals, in terms of consecutive values of
the variable $k$:

\begin{center}
\textbf{8, 9, 9, 9, 8, 9, 9, 9, 8, 9, 9 } for a sum of 96 values of the
variable $k$.
\end{center}

As we can see, the first interval of 8 values has been neglected: in it the
function $\Omega(I_1 )$ has value zero, but such interval has a
peculiarity, as exactly at interval 31 the above mentioned \textit{series
  }is interrupted and an interval of \textbf{8} values takes place
\footnote{ For precision it is the last 9, belonging to the second group of
three consecutive 9 , that is of the eighth term of the \textit{series.}}
of one of \textbf{9, }giving origin again to the same sequence: the
following 11 values of such function are determined by the extent of
intervals, typical of the \textit{series.}

Now this event is regularly repeated and precisely every \textbf{40} and
\textbf{51} intervals: for instance the first substitutions happen at the
intervals: 31, 71, 122, 162, 213, 253, 293, 344, 384, 435, 475,
526,566,617,657,697,748,788,839,879,930,970,1021,1061,1112,1152 and, all of
them, in advantage of the 8 value intervals.

The first real effect of this phenomenon on a procedure that allows
to calculate exactly the last value of $k$, present in a determined
interval, foreseeing a constant and complete repetition of the
\textit{numerical series} of 11 intervals, is exactly quantifiable
at the interval 122, where it is practically possible to verify that
the last value of $k$ inside it, results lower of a \textbf{unity}:
this means that we could obtain the same result if at the interval
122, excluding the initial one, there were \textbf{120} integer
intervals and one reduced of a unity, because of one only
interruption.

Growing $k$, the substitutions immediately determinate a further but stable
effect, for instance, at the interval 617 it's possible to verify that the
last value of $k$, present in the interval is of \textbf{two }unities inferior
than the computable one and, going further it's at the interval 1112, that
it is possible to verify that the last value of $k$, present in the
interval, is of \textbf{three} unities inferior than the computable one,
and nothing short, for many following intervals it's possible to observe
that such reduction of a unity is noticeable exactly every 495 intervals
for indeed 10 times (\textbf{first case}) and every 484 intervals only once
(\textbf{second case}).

In the \textbf{first case} we can get the same results by hypothesizing the
constant presence of 494  integral intervals and one, reduced of a
unity and this because \textbf{11} substitutions produce inside 495
intervals \textbf{38 }complete \textit{series} of 96 values of $k$ and
\textbf{11 }\textit{series, }reduced to only 61 values of $k$; under these
conditions it's immediate to verify that we obtain 4319 values of $k$, that
are inferior of one only unity in relation with the possible values (4320)
inside 495 intervals, in case there were exclusively 45 numeric complete
\textit{series}.

In the \textbf{second case} we can get the same results by assuming the
constant presence of 483  integral intervals and one, reduced of a
unity and this because \textbf{11} substitutions produce inside 484
intervals \textbf{37 }complete \textit{series} of 96 values of $k$ and
\textbf{11 }\textit{series, }reduced to only 61 values of $k$; under these
conditions it's immediate to verify that we obtain 4223 values of $k$, that
are inferior of one only unity in relation with the possible values (4224)
inside 484 intervals, in case there were exclusively 44 numeric complete
\textit{series.}
\newpage
The peculiarity of the first interval of 8 values, for which the function
$\Omega(I_\ell )$ is worthless, consists therefore in the fact that, as
it was determined before the birth of the sequence, characteristic of the
\textit{series}, it is possible to think that at the origin, the number of
integer intervals corresponds to \textbf{1.}

We specify that such characteristics, even if relevant, are minor, both in
relation to the continuous and repetitive presence of the above described
\textit{series} and in relation with the largeness of the single intervals,
which never descend under \textbf{8} values\footnote{ This means that every
interval $I_\ell $ will never be empty.} for very great values of $k$.

The ratio of the number of integers inside the numeric \textit{series }
with the 11 intervals corresponds to $8,\overline {72} $ and this value
diminishes in a little meaningful way while $k$ tends to infinity, if we
consider the ratio of the value of a very great $k$ with its own belonging
interval $I_\ell $: a sufficiently precise value is obtainable with the
following expression
$\frac{\mbox{96}}{\mbox{11}}-\frac{\mbox{1}}{\mbox{494}}\;$, where we deduce
by the denominator of the second fraction the importance of the number
\textbf{494}.

Successively we give (Fig. 3) the graphs, related both to the
family of <^s> functions, or better, to the set of the arcs
belonging to them, and to the auxiliary function $\Omega(I_\ell )$
and successively in particular the development of the graph of a
couple of functions, characteristic of a precise interval $I_\ell
$ (Fig. 4) and of the pointers than put into evidence the
behaviour, growing $k$.

This last behaviour results particularly evident, by examining the
various displays, produced by \textit{DERIVE}, in relation with the
first 11 intervals, typical of the \textit{standard series}, while
for the greater values of $k$ it's necessary to outdistance in an
appropriate way the intervals on which to do a comparison to have a
further confirmation of such behaviour.

Now if $k$ tends to infinity it's possible to compute the limit
towards which the generic <^s> function, that in this case
represents the fusion, to infinity, of two arcs of
\textit{following} <^s> functions, that refer themselves to a
hypothetical and extreme interval $I_\ell $.

\medskip
\centerline{\epsfig{figure=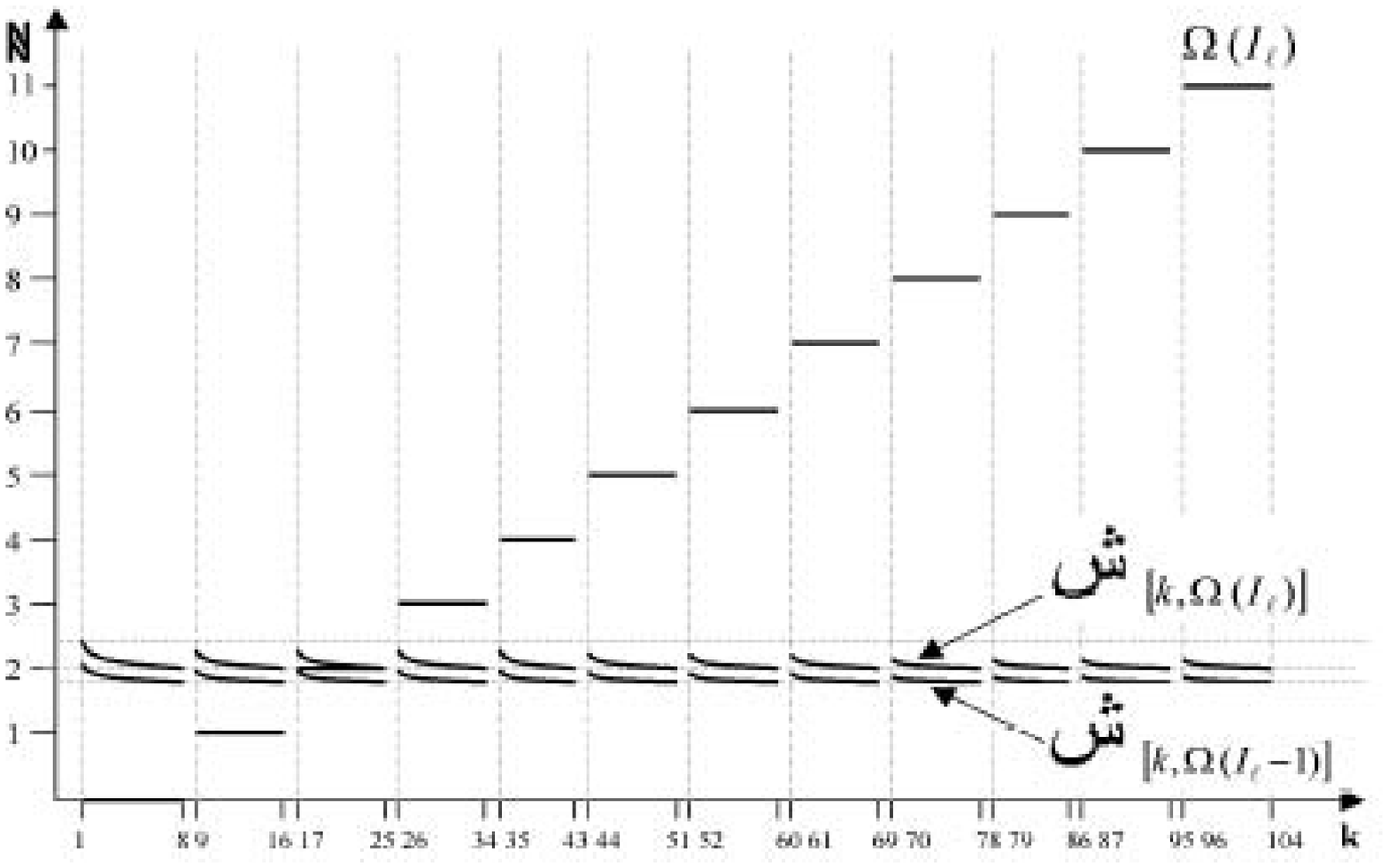,width=11cm,height=6.75cm}}
\medskip
\baselineskip=0.20in
\begin{center}
{\bf Fig. 3}
\end{center}

Keeping in mind what noticed in the first and second case,
previously described we can, first of all, calculate the two
following values, by ``excess'' and ``defect'' of the ratio of a
value of $k$, with the value of the integer function $\Omega(I_\ell )$,
determined as a function of the $I_\ell $
interval, containing the same value of $k$:

\baselineskip=0.20in

\[ 1^{st}\mbox{case\textbf{:}}\frac{k}{\Omega(I_\ell)}\prec
\frac{4319}{495}=8,7\overline{25}\]
\[ 2^{nd}\mbox{case\textbf{:}}
\frac{k}{\Omega(I_\ell)}\succ\frac{4223}{484}=8,72\overline
{52066115702479333884297} \]

These values, by considering the frequency of the two cases ( 10
times the first case and once the second case), permit to the
approximate the ratio with:
\[
\frac{k}{\Omega(I_\ell)}\cong \frac{10\ast
8,7\overline{25}+8,72\overline{52066115702479333884297} }{11}=
8,7252483512814091326488=
\]
\[
=\frac{96}{11}-\frac{1}{493,9793814432989690721649}\cong
\frac{96}{11}-\frac{1}{494}
\]

\baselineskip=0.20in

In this way, according to the analysis of the progress of the
various and \textit{following} <^s> functions, we can furnish only
an esteem of the ratio $\frac{k}{\Omega \;(I_\ell )\;}$.

If we want to determinate the effective value of such ratio, it's
necessary to compare the values of the two quantities for
sufficiently great values of $k$.

Keeping in mind the algorithm described at the end of paragraph 1,
which puts in evidence the continuous oscillation of the rational
values of the \textit{following} <^s> functions, in proximity of
the integer number 2, we will compare the values of the ratio
$\frac{k}{\Omega \;(I_\ell )\;}$ as $k$ grows with the value:
$\left( {3-2 \mathord{\left/ {\vphantom {2 {\log 2}}} \right.
\kern-\nulldelimiterspace} {\log 2}} \right)^{-1}$.

\medskip
\centerline{\epsfig{figure=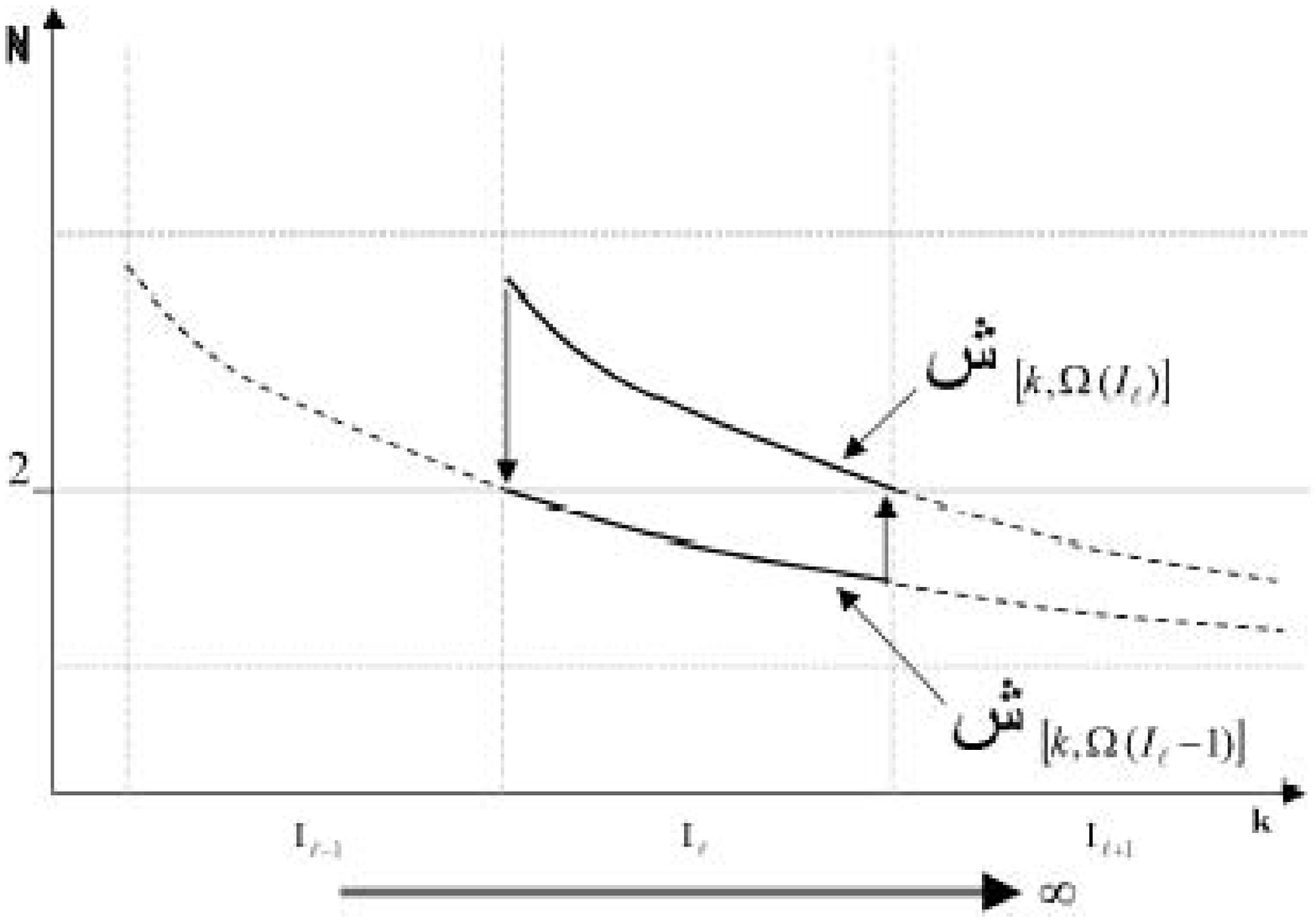,width=10cm,height=7cm}}
\medskip
\baselineskip=0.20in
\begin{center}
{\bf Fig. 4}
\end{center}

More precisely, considering the values raised to the $18^{th}$power ($k_1
=10^{18}$, that is one quintillion), raised to the $33^{rd}$ power ($k_2
=10^{33}$, that is one decillion) and indeed raised to the $63^{rd}$ power
($k_3 =10^{63}$, that is one vigintillion) it's possible to observe that the
real value of the ratio $\Psi \,(k)=\frac{k}{\Omega \;(I_\ell )\;}$
tends to the previously shown value; in fact in the three various cases it's
possible to determinate in order what follows:
\[
\Rightarrow\Psi(k_1)=\frac{k_1}{\Omega(I_\ell)}
-\left({3-2\mathord{\left/{\vphantom{2{\log2}}}\right.
\kern-\nulldelimiterspace}{\log 2}}\right)^{-1}\prec5\cdot10^{-19};
\Rightarrow\Psi(k_2)-\left({3-2\mathord{\left/{\vphantom{2{\log2}}}\right.
\kern-\nulldelimiterspace}{\log
2}}\right)^{-1}\prec5\cdot10^{-34};\]
\[\Rightarrow\Psi(k_3)-\left({3-2\mathord{\left/{\vphantom{2{\log2}}}\right.
\kern-\nulldelimiterspace}{\log 2}}\right)^{-1}\prec5\cdot10^{-64}
\]

Definitively, going to the limit for $k$, which tends to infinity
we can expect that:
\[
\lim _{k\to \infty } \left( {\frac{k}{\Omega \;(I_\ell )}}
\right)=\left( {3-2 \mathord{\left/ {\vphantom {2 {\log 2}}}
\right. \kern-\nulldelimiterspace} {\log 2}} \right)^{-1}
\]

This result has obviously an immediate consequence in the calculation of the
following limit:

\[\qquad\quad\lim _{k\to \infty }  \quad \mbox{<^s>} \left[ {k,\Omega
\;(I_\ell )} \right]\;=\lim _{k\to \infty } \;\left(
{1+\frac{1}{3k-\Omega \;(I_\ell )}\mbox{ }} \right)^{2k+1}\;=\]

\[\quad = \lim _{k\to \infty }
e^{\frac{2k+1}{3k-\Omega \;(I_\ell )}}\;\; = e^{\lim_{k\to\infty}
\frac{2k+1}{3k-k\left({3-\frac{2}{\log \,2}} \right)}}\;\;=e^{\log
\,2}\;\;=\;\;2\]

\vspace{3mm}

In conclusion, it happened that the natural number 2 (by the integer
constant function Y=2) can represent the separation element of two
arcs, belonging to two \textit{following} <^s> functions, which
therefore can be identified with two contiguous classes\footnote{
The contiguous classes
  are meant represented by two groups of rational non-integer numbers,
  greater and smaller than 2, separated by the rational integer number 2,
which obviously doesn't belong to any of the two classes.}, which tend to
approach indefinitely.

Even if we observe that the validity of the \textit{fundamental }
theorem is actually included in the construction algorithm\footnote{
Let's observe that the algorithm determinates for the generic <^s>
function an inferior limit: it's established an inferior extreme
which corresponds to the integer number 2.} of the family of <^s>
 functions, keeping obviously in mind the characteristic of
monotonicity (see previous paragraph and the following) of the
generic <^s>  function, by a more appropriate notation, due to
Iverson, further we will give the elements to get a strict proof of
the theorem in the modern sense of the term.

Let's extend the dependence of the integer step function $\Omega
\;(I_l )$ to the real field and let's use the following
definition:

\vspace{3mm}

\[ \Omega\left( x\right)=\mbox{min}\left\{
{k\in}\right.\mathbf{N}\mbox{: S}_{k+1}\left( x\right)\ge
2\left.\right\}\mbox{; }x\in\mathbf{R^{+}}\mbox{ and where S}_k
\left( x\right)=\left({1+\frac{1}{3x-k+1}} \right)^{2x+1}\]

\vskip 0.1truecm

\newpage

By simple algebraic passages we have that:
\begin{equation}
\label{eq2} \Omega\;\left( x \right)=\left\lceil
{3x-\frac{1}{2^{\frac{1}{2x+1}}-1}} \right\rceil
\end{equation}
where $\left\lceil {\,x} \right\rceil $ means the smallest
integer, greater than $x$ or equal to it.

So, from such a particular definition of $\Omega\;\left(
x \right)$ the \textit{fundamental }theorem immediately
derives.

\vskip 0.1truecm

\vspace{3mm}

In fact the <^s>  function possesses the following explicit formula:

\[ \mbox{S}\left( x \right)=\left(
{1+\frac{1}{3x-\Omega \left( x \right)}} \right)^{2x+1}=\left(
{1+\frac{1}{3x-\left\lceil {3x-\frac{1}{2^{\frac{1}{2x+1}}-1}}
\right\rceil }} \right)^{2x+1}\]\

\noindent and passing to the sequence we have:
\[\mbox{S}\left( n \right)\;\,=\;\; \quad \left(
{\mbox{1}+\frac{\mbox{1}}{\left\lfloor
{\frac{1}{2^{\frac{1}{2n+1}}-1}} \right\rfloor }} \right)^{2n+1}\]

\noindent where $\left\lfloor x \right\rfloor $ shows the greatest
integer smaller or equal to $x$.

From now on we will use the locution of `` special function <^s> ''
for the discrete function $\mbox{<^s>} \left[ {k,\Omega \;(I_\ell
)} \right]\;$, defined by the formula (\ref{eq1}).

Extending the field of definition of the variable $k$ to the real
positive numbers, it's possible to notice that such function,
being represented by the union of continuous arcs (all above the
straight line of height 2 ) is actually assimilable to a piecewise
continuous function.

\vspace{8mm}

\setnashbf
\begin{center}
3. AN APPROACH OF THE SPECIAL FUNCTION <^s> WITH THE
EULERIAN\footnote{ Chapter II, [5].} GAMMA FUNCTION $ \mathbf{\Gamma} $
\end{center}
\setnash

\vspace{6mm}

 One of the fundamental characteristics of the
Eulerian \textit{gamma} function is a certain condition of
monotonicity, which is the fact that the function, which can be
intended as the most spontaneous extension of the factorial $n\,!
=1\cdot 2\cdot 3\cdot \cdot \cdot n$, out the field of the natural
numbers, is logarithmically convex.

\newpage

If we consider the expression (\ref{eq1}) and we calculate its first derivative,
we have:

\[ \frac{d}{dk} \mbox{<^s>} \left[ {k,\Omega \;(I_\ell )} \right]\;=\]
\begin{equation}
\label{eq3} \mbox{<^s>} \left[ {k,\Omega \;(I_\ell )} \right]\;
\cdot \left[ {2\log \left({1+\frac{1}{\left( {3k-\Omega \;(I_\ell )}
\right)}}\right){\kern 1pt}{\kern 1pt}-\frac{3\cdot \left( {2k+1}
\right)} {\left({3k-\Omega \;(I_\ell )} \right)\left( {3k-\Omega
\;(I_\ell )+1} \right)}}\right]
\end{equation}
$\forall \;\;k\in \mathbf{N}$ and $k>\frac{\Omega \;(I_\ell )}{3}$ it is
always negative, therefore the function results monotonic and decreasing,
and besides it's interesting to observe that, growing the $k$ value, its
absolute value diminishes.

That being stated we can also verify that the special function <^s>
has the same characteristic of monotonicity of the Eulerian function
$\Gamma $: it's sufficient in fact to verify that the second
derivative of the logarithm of the same <^s> function is positive.

Such result is immediate, in fact for the (\ref{eq3}), being:

\[ \frac{d^2}{dk^2}\log  \mbox{<^s>} \left[ {k,\Omega \;(I_\ell )} \right]\;=
\frac{d}{dk}\left\langle \right. \quad \frac{d}{dk} \mbox{<^s>}
\left[ {k,\Omega \;(I_\ell )} \right]\; \raise0.7ex\hbox{}
\!\mathord{\left/ {\vphantom {
}}\right.\kern-\nulldelimiterspace}\!\lower0.7ex\hbox{}
\;\mbox{<^s>} \left[ {k,\Omega \;(I_\ell )} \right]\; \quad \left.
\right\rangle \]

\noindent we have that:

\[ \frac{d^2}{dk^2}\log\mbox{<^s>}\left[{k,\Omega(I_\ell)}
\right]=\]
\[
\frac{d}{dk}\left[{2\log{\kern 1pt}\left({1+\frac{1}{\left(
{3k-\Omega(I_\ell)}\right)}}\right) -\frac{3\cdot\left({2k+1}
\right)}{\left({3k-\Omega(I_\ell)}\right)\cdot\left({3k-\Omega(I_\ell)+1}
\right)}}\right]=\]

\[ = \left[ {\frac{12\cdot \left( {3k\cdot \Omega \;(I_\ell )-\Omega \;(I_\ell
)^2} \right)+6\cdot \left( {6k-\Omega \;(I_\ell )} \right)+9}{\left(
{3k-\Omega \;(I_\ell )} \right)^2\cdot \left( {3k-\Omega \;(I_\ell )+1}
\right)^2}} \right]\]

Since, for definition we have $k\succ \Omega \;(I_\ell )\;$,$\forall \;k\in
\mathbf{N}$, we have in conclusion:

\[ \frac{d^2}{dk^2}\log \mbox{<^s>} \left[ {k,\Omega \;(I_\ell )}
\right]\;\;\;\succ \;0\]

This demonstrates that also the special function <^s> is
logarithmically convex.

The special function <^s> is actually an exponential general
function; that being stated, considering its base, defined by the
analysis of its behaviour to infinity, we consider the following two
particular expressions:

\[ \prod\limits_{k=1}^n
\left( {1+\frac{\log \,\,2}{2k}} \right) \;\;\;\mbox{and}\;\;\;
\sum\limits_{k=1}^n \left( {\frac{\log \,\,2}{2k+\log 2}}
\right)\]

In the previous paragraph we have practically put into evidence that
for great values of $k$ the rational non-integer term of the base of
the special function <^s> tends ($\sim )$ to:

\begin{equation}
\label{eq3_1}
\frac{1}{3k-\Omega \;(I_\ell )} \sim  \frac{\log
  \,\,2}{2k}
\end{equation}

\noindent consequently, as there is the following relation $\Omega
\;(I_\ell -1)=\Omega \;(I_\ell )-1$, we have:

\begin{equation}
\label{eq3_2}
\frac{1}{3k-\Omega \;(I_\ell -1)} \sim  \frac{\log
  \,\,2}{2k+\log 2}
\end{equation}

The value (\ref{eq3_1}), if fixed for the base of the generic <^s>
function, is very well fit to ``interpolate'' in the continuous
field the union of the various arcs of the family of <^s> functions,
got in the discrete field, above the straight line of height 2.

By this last value it's possible to calculate and verify\footnote{ $\Gamma
\left( {x+n} \right)=x\cdot \left( {x+1} \right)\cdot \left( {x+2}
\right)\cdot \cdot \cdot \left( {x+n-1} \right)\cdot \Gamma \left( x
\right)$ with $n\in \mathbf{N}$ and $x\succ 0$} the following sizeable
expression:

\begin{equation}
\label{eq4}
\prod\limits_{k=1}^n
\left( {1+\frac{\log 2}{2k}} \right)
\quad
=\frac{\Gamma \left( {n+1+\frac{\log \,2}{2}} \right)}{\Gamma \left( {n+1}
\right)\cdot \Gamma \left( {1+\frac{\log \,2}{2}} \right)}
\end{equation}

The value (\ref{eq3_2}), on the contrary, if fixed for the base of
the generic <^s> function, is very well fit to ``interpolate'' in
the continuous field the union of the various arcs of the family of
<^s> functions, got in the discrete field, below the straight line
of height 2.

By it, it's possible to define the following partial sum:
\begin{equation}
\label{eq5}
\sum\limits_{k=\mbox{1}}^n
\frac{\log \,\,2}{2k+\log 2}\quad =
\quad
\frac{\log \,\,2}{2}
\sum\limits_{k=1}^n
\frac{1}{k+\frac{\log 2}{2}}
\end{equation}

From here, keeping in mind the formulae of recurrence of the logarithmic
derivative of the \textit{gamma, }function named \textit{digamma} :

\[ \psi(x+1)=\psi(x)+\frac{1}{x} ;
\psi(x+n)=\psi(x)+\frac{1}{x}+\frac{1}{x+1}+\frac{1}{x+2}+.....+
\frac{1}{x+n-1}\mbox{ with }n\ge 1\]

\noindent we have, giving x the value $1+\frac{\log \,2}{2}$, the
following identity:
\[
\sum\limits_{k=\mbox{1}}^n
\frac{\log \,\,2}{2k+\log 2}
\quad
=\frac{\log \;2}{2}\left[ {\psi \;\left({n+1+\frac{\log \;2}{2}}\right)-\psi
\;\left({1+\frac{\log \;2}{2}}\right)} \right]
\]

Now, going to the limit for $n\to \infty $, for the known properties of the
\textit{digamma} function, we have that the following series is divergent,
that is:
\begin{equation}
\label{eq6}
\sum\limits_{k=1}^\infty
\frac{\log \,\,2}{2k+\log 2}\quad =
\quad
\infty
\end{equation}

That being stated, the importance of the result (\ref{eq4}) must
be evaluated above all according to the following identity:
\[
\mbox{ }\left( {1+\frac{\log 2}{2k}} \right)\quad =
\quad
\left( {\mbox{1}-\frac{\log \,\,2}{2k+\log 2}} \right)^{-1}
\]

From this last identity we can deduce:
\begin{equation}
\label{eq7}
\prod\limits_{k=1}^n
\left( {1-\frac{\log \,\,2}{2k+\log 2}} \right)\;\;=\
\frac{\Gamma \left( {n+1} \right)\cdot \Gamma \left( {1+\frac{\log \,2}{2}}
\right)}{\Gamma \left( {n+1+\frac{\log \,2}{2}} \right)}
\end{equation}
and considering the infinite product:

\[ \prod\limits_{k=\mbox{1}}^\infty
\left( {\mbox{1}-\frac{\log \,\,2}{2k+\log 2}} \right)\;\;\]

\noindent diverges \footnote{ Chapter II, pag. 33, [6].} to zero.

Such conclusion can be immediately verified also by applying the result
(\ref{eq6}), modified in the sign and considering the following:

\vspace{3mm}

\textbf{Theorem}: \textit{Supposed that }$-1\prec a_n \le 0$\textit{ and
  the series }$\sum {a_n } $\textit{ is divergent,
  then the infinite product }$\prod {\left( {1+a_n } \right)} $
\textit{ diverges to zero.}

\begin{center}
  \textit{Proof}
 \end{center}

\quad Assumed $b_n =-\,a_n \Rightarrow \quad \;0\le \;\;b_n \prec
\,\;1$, as for $0\le x\prec 1$ results $1-x\le e^{\mbox{-}x}$, we
can write:

\begin{equation}
\label{eq8}
0\prec \quad P_n
=\prod\limits_{r=1}^n {\left(
{1+a_r } \right)\quad } \;\le \;\quad
e^{-\left( {b_1 +b_2
+.....+b_n } \right)}
\end{equation}
then if the series $\sum {a_n }$ is not convergent and necessarily
diverges to $-\infty $, from (\ref{eq8}) we have that $P_n \to 0$
that is the infinite product diverges to zero.

Besides, keeping in mind that, supposed $a$ and $b$ non negative, the following
relation is valid:
\[ \prod\limits_{r=\mbox{1}}^\infty
\frac{r\cdot \left( {r+a+b} \right)}{\left(
{r+a} \right)\cdot \left( {r+b} \right)}=\frac{\Gamma
\left( {1+a} \right)\cdot \Gamma \left( {1+b} \right)}{\Gamma \left( {1+a+b}
\right)}\]

\noindent verifiable, applying the known formula by Euler:
\[
\Gamma \left( x \right)=\lim _{n\to \infty } \frac{n^x\cdot n!}{x\cdot
\left( {x+1} \right)\cdot \left( {x+2} \right)\cdot \cdot \cdot \cdot \left(
{x+n} \right)}
\]
and considering the limit for $n\to \infty $ of the development of the
following finished product:
\[
\prod\limits_{r=\mbox{1}}^{n+1}
\frac{r\cdot \left( {r+a+b} \right)}{\left(
{r+a} \right)\cdot \left( {r+b} \right)}
=\frac{\mbox{(1}+a+b)\cdot \mbox{(2}+a+b)\cdot
\cdot \cdot \left( {\mbox{1}+a+b+n}
\right)}{(n^{1+a+b})\cdot n!}\cdot\]
\[ \cdot\frac{(n^{1+a})\cdot
n!}{\mbox{(1}+a)\cdot (2+a)\cdot \cdot \cdot \left(
{\mbox{1}+a+n} \right)}
\cdot \frac{(n^{1+b})\cdot n!}{\mbox{(1}+b\mbox{)}\cdot (2+b)\cdot \cdot
\cdot \left( {\mbox{1}+b+n} \right)}\cdot \frac{n+1}{n}
\]
we can state that:
\[ \frac{\Gamma \left( {1+n} \right)\cdot \Gamma \left(
{1+\frac{\log 2}{2}} \right)}{\Gamma \left( {1+n+\frac{\log 2}{2}}
\right)}= \quad \prod\limits_{r=\mbox{1}}^\infty \frac{r\cdot
\left( {r+n+\frac{\log 2}{2}} \right)}{\left( {r+n}
\right)\cdot \left( {r+\frac{\log 2}{2}} \right)}\]

By exploiting the properties of the \textit{gamma} function, or resorting
to the definition of the famous hypergeometric series, it's also possible
to verify the following:
\begin{equation}
\label{eq9} \frac{\Gamma \left( {1+n} \right)\cdot \Gamma \left(
{1+\frac{\log 2}{2}} \right)}{\Gamma \left( {1+n+\frac{\log 2}{2}}
\right)}= \sum\limits_{r=0}^\infty \frac{\left( {-1}
\right)^r\cdot \left( {\frac{\log 2}{2}} \right)\cdot \left(
{\frac{\log 2}{2}-1} \right)\cdot \cdot \cdot \left( {\frac{\log
2}{2}-r} \right)}{\left( {1+n+r} \right)\cdot r\;!}
\end{equation}

In fact, adopting the classic symbolism for the hypergeometric
function\footnote{ Chapter III, [5].} :
\[
F\left( {a,b;c;x} \right)=1+\frac{a\cdot b}{c\cdot 1\;!}\cdot x+\frac{a\cdot
\left( {a+1} \right)\cdot b\cdot \left( {b+1} \right)}{c\cdot \left( {c+1}
\right)\cdot 2\;!}\cdot x^2+\]
\[ +\frac{a\cdot \left( {a+1} \right)\cdot \left(
{a+2} \right)\cdot b\cdot \left( {b+1} \right)\cdot \left( {b+2}
\right)}{\;c\cdot \left( {c+1} \right)\cdot \left( {c+2} \right)\cdot
3\;!}\cdot x^3+\;....
\]
we can exploit the possibility to express it in terms of \textit{gamma}
functions, considering the hypergeometric integral, that is:
\[
F\left( {a,b;c;x} \right)=\frac{\Gamma (c)}{\Gamma \left( b \right)\cdot
\Gamma \left( {c-b} \right)}\cdot \int\limits_0^1 {t^{b-1}\cdot } \left(
{1-t} \right)^{c-b-1}\cdot \left( {1-t\cdot x} \right)^{-a}dx
\quad
(\Re c\succ \Re b\succ 0)
\]

From this, considering the limit for $x\to 1^-$ (Abel's theorem) and
exploiting the properties of the Eulerian \textit{Beta} function, we obtain
the important relation of the hypergeometric Gauss's series:
\[
F\left( {a,b;c;\;1} \right)=\frac{\Gamma (c)\cdot \Gamma \left( {c-a-b}
\right)}{\Gamma \left( {c-a} \right)\cdot \Gamma \left( {c-b} \right)}
\quad
(c\ne 0,-1,-2,...,\;\;\Re \left( {c-a-b} \right)\succ 0)
\]

\vspace{2mm}

 With the positions $a=-\frac{\log 2}{2}$, $b=n$ and $c=n+1$,
as all the required limitations are satisfied for the parameters, we
easily reach the (\ref{eq9}).

In conclusion, starting from the characteristic base of the generic
<^s> function, defined by the analysis of its behaviour to infinity,
we have also stated the following sizeable relation between an
infinite product and a numeric series:
\[
F\left( {-\frac{\log 2}{2},\,\,n;\,\;n+1;\;\;1} \right)=
\]
\[ \prod\limits_{r=\mbox{1}}^\infty \frac{r\cdot \left(
{r+n+\frac{\log 2}{2}} \right)}{\left( {r+n} \right)\left(
{r+\frac{\log 2}{2}} \right)}=\sum\limits_{r=0}^\infty \frac{\left( {-1}
\right)^r\cdot \left( {\frac{\log 2}{2}} \right)\cdot \left( {\frac{\log
2}{2}-1} \right)\cdot \cdot \cdot \left( {\frac{\log 2}{2}-r}
\right)}{\;\left( {1+n+r} \right)\,{\kern 1pt}\cdot {\kern 1pt}r\;!}\]

\vspace{3mm}

In the real field such result is graphically represented by Fig. 5,
shown below, where it is evident that the $x$ axis represents a
horizontal asymptote for $x\to +\infty $.

\newpage

Besides, in the field of the non-negative real numbers, both the
infinite product and the series are absolutely and uniformly
convergent.

\medskip
\centerline{\epsfig{figure=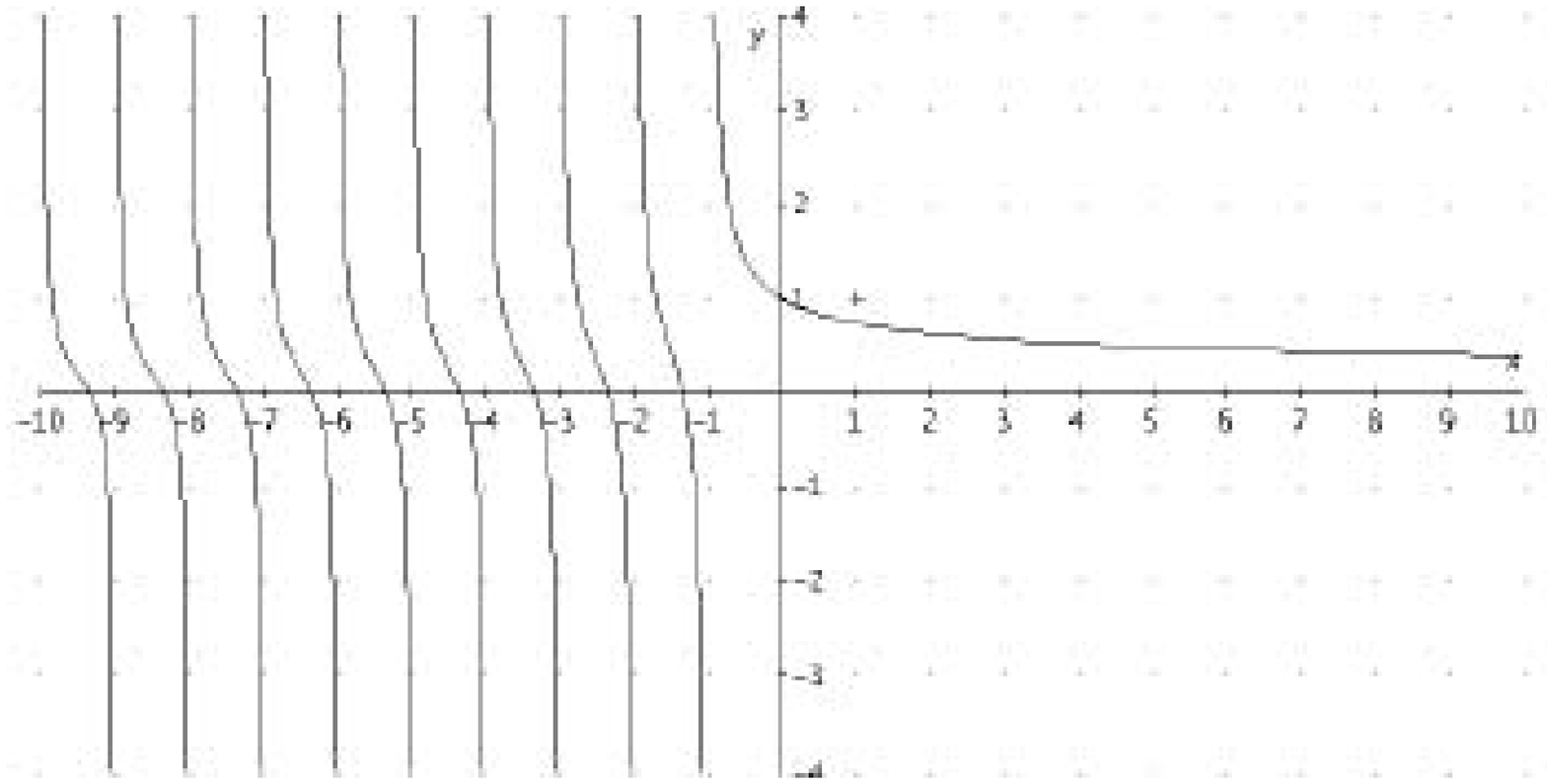,width=14cm,height=7cm}}
\medskip
\baselineskip=0.20in
\begin{center}
{\bf Fig. 5}
\end{center}

\vspace{8mm}
\setnashbf
\begin{center}
4. THE SPECIAL FUNCTION <^s> AND THE SECOND-ORDER EULERIAN
NUMBERS\footnote{ pag. 247-251,[3].}
\end{center}

\setnash

In the second paragraph we have shown that one of the most important
characteristics in the construction of the family of <^s> functions
is related with the repetitive presence of 11 specific intervals
$I_\ell $ of the integer variable $k$: we have indeed seen that the
extension of such 11 intervals is always characterized by the
following numeric \textit{series}:

\begin{center}
\textbf{8, 9, 9, 9, 8, 9, 9, 9, 8, 9, 9 } for a sum of 96 values of the
variable $k$.
\end{center}

But it also happened that such \textit{series} undergoes inside a precise
number of integer intervals some interruptions, as 8 value intervals
substitute some 9.

More precisely this event is regularly repeated every \textbf{40} and
\textbf{51} intervals and it's possible to observe that such intervals are
aggregable so that to give origin to two groups: the first one formed by
\textbf{484 }intervals, produced by \textbf{7} interruptions at a distance
of \textbf{40} intervals and \textbf{4 }interruptions at a distance of
\textbf{51 }intervals, the second one formed by \textbf{495 }intervals,
produced by \textbf{6 }interruptions at a distance of \textbf{40} intervals
and \textbf{5} interruptions at a distance of \textbf{51 }intervals.

For each group of \textbf{484} intervals we besides note always \textbf{10}
groups of \textbf{495 }intervals: this corresponds on average to a number of
\textbf{494 }intervals for all the groups.

In conclusion the characteristics of such phenomenon, keeping also
in mind what explained in the initial part of the second
paragraph, allow correctly some estimates for the integer values
of $k$, referred to an interval $I_\ell $, if we suppose the
effective presence, even if virtual, of the following sequence of
integer intervals:
\vspace{3mm}

\textbf{1,120,494,494,{\ldots},494, }$\approx
$\textbf{494}\footnote{ $\approx $\textbf{494:} with such notation
we want to put into evidence that for great values of $k$ the
estimable number of integer intervals is, in some very near case,
to such value, in fact it's as if every 495.000 intervals, two
others of them, on average, undergo a reduction of one unity (from
9 to 8), but on the whole this further phenomenon is absolutely
neglectable.}\textbf{,494,..,494,{\ldots},}\textbf{494,..,494,{\ldots},}
$\approx $\textbf{494} (ad infinitum).

\vspace{3mm}

This interpretation is besides confirmed by the check carried out
on the limit of the ratio $\raise0.7ex\hbox{$k$} \!\mathord{\left/
{\vphantom {k {\Omega \;(I_\ell
)}}}\right.\kern-\nulldelimiterspace}\!\lower0.7ex\hbox{${\Omega
\;(I_\ell )}$}$.

Now, as these numbers are closely connected to the nature of the
complete and periodic sequence of the numerical series shown above,
we have supposed that the numbers \textbf{1 , 120 }and\textbf{ 494}
can belong to a category of \textit{special numbers;} particularly a
research pointed to verify this hypothesis, has implied the
following curious discovery even if it is always a conjecture on an
almost light and remote connection: the numbers 1, 120 and 494
belong to the family of the so called ``second-order Eulerian
numbers'', which result important for the tight connection that they
have with Stirling's numbers.

They satisfy similar recurrence to the characteristic one of the ``ordinary
Eulerian numbers'' , which are useful, above all because they give a
connection among ordinary powers and consecutive binomial coefficients.

To be clearer we show below the recurrence, which characterizes the
second-order Eulerian numbers, showing how it's possible, by it, to
produce the three special numbers, typical of the special function
<^s>.

\vspace{6mm}

In fact we have:

\vspace{6mm}

\qquad\qquad\qquad\quad $\left\langle {\left\langle
{\begin{array}{l}
 n \\
 k \\
 \end{array}} \right\rangle } \right\rangle =(k+1)\left\langle
{\left\langle {\begin{array}{l}
 n-\mbox{1} \\
 \;\,\,k \\
 \end{array}} \right\rangle } \right\rangle +(2n-1-k)\;\left\langle
{\left\langle {\begin{array}{l}
 n-\mbox{1} \\
 k-\mbox{1} \\
 \end{array}} \right\rangle } \right\rangle $ ;

\vspace{6mm}

\qquad\qquad\qquad\quad $\left\langle {\left\langle
{\begin{array}{l}
 n \\
 0 \\
 \end{array}} \right\rangle } \right\rangle =\textbf{1}
\quad \forall \mbox{n}\ne 0\quad ; \quad \left\langle
{\left\langle {\begin{array}{l}
 n \\
 n \\
 \end{array}} \right\rangle } \right\rangle =0$ for $n\ne 0$

\newpage
\[
\left\langle {\left\langle {\begin{array}{l}
 2 \\
 1 \\
 \end{array}} \right\rangle } \right\rangle =(2)\left\langle {\left\langle
{\begin{array}{l}
 1 \\
 1 \\
 \end{array}} \right\rangle } \right\rangle +\left( 2 \right)\left\langle
{\left\langle {\begin{array}{l}
 1 \\
 0 \\
 \end{array}} \right\rangle } \right\rangle =2
\]
\[
\left\langle {\left\langle {\begin{array}{l}
 3 \\
 2 \\
 \end{array}} \right\rangle } \right\rangle =(3)\left\langle {\left\langle
{\begin{array}{l}
 2 \\
 2 \\
 \end{array}} \right\rangle } \right\rangle +\left( 3 \right)\left\langle
{\left\langle {\begin{array}{l}
 2 \\
 1 \\
 \end{array}} \right\rangle } \right\rangle =6
\]
\[
\left\langle {\left\langle {\begin{array}{l}
 4 \\
 3 \\
 \end{array}} \right\rangle } \right\rangle =(4)\left\langle {\left\langle
{\begin{array}{l}
 3 \\
 3 \\
 \end{array}} \right\rangle } \right\rangle +\left( 4 \right)\left\langle
{\left\langle {\begin{array}{l}
 3 \\
 2 \\
 \end{array}} \right\rangle } \right\rangle =24
\]
\[
\left\langle {\left\langle {\begin{array}{l}
 5 \\
 4 \\
 \end{array}} \right\rangle } \right\rangle =(5)\left\langle {\left\langle
{\begin{array}{l}
 4 \\
 4 \\
 \end{array}} \right\rangle } \right\rangle +\left( 5 \right)\left\langle
{\left\langle {\begin{array}{l}
 4 \\
 3 \\
 \end{array}} \right\rangle } \right\rangle =\textbf{120}
\]
\vspace{6mm}
\[
\left\langle {\left\langle {\begin{array}{l}
 2 \\
 1 \\
 \end{array}} \right\rangle } \right\rangle =(2)\left\langle {\left\langle
{\begin{array}{l}
 1 \\
 1 \\
 \end{array}} \right\rangle } \right\rangle +\left( 2 \right)\left\langle
{\left\langle {\begin{array}{l}
 1 \\
 0 \\
 \end{array}} \right\rangle } \right\rangle =2
\]
\[
\left\langle {\left\langle {\begin{array}{l}
 3 \\
 1 \\
 \end{array}} \right\rangle } \right\rangle =(2)\left\langle {\left\langle
{\begin{array}{l}
 2 \\
 1 \\
 \end{array}} \right\rangle } \right\rangle +\left( 4 \right)\left\langle
{\left\langle {\begin{array}{l}
 2 \\
 0 \\
 \end{array}} \right\rangle } \right\rangle =8
\]
\[
\left\langle {\left\langle {\begin{array}{l}
 4 \\
 1 \\
 \end{array}} \right\rangle } \right\rangle =(2)\left\langle {\left\langle
{\begin{array}{l}
 3 \\
 1 \\
 \end{array}} \right\rangle } \right\rangle +\left( 6 \right)\left\langle
{\left\langle {\begin{array}{l}
 3 \\
 0 \\
 \end{array}} \right\rangle } \right\rangle =22
\]
\[
\left\langle {\left\langle {\begin{array}{l}
 5 \\
 1 \\
 \end{array}} \right\rangle } \right\rangle =(2)\left\langle {\left\langle
{\begin{array}{l}
 4 \\
 1 \\
 \end{array}} \right\rangle } \right\rangle +\left( 8 \right)\left\langle
{\left\langle {\begin{array}{l}
 4 \\
 0 \\
 \end{array}} \right\rangle } \right\rangle =52
\]
\[
\left\langle {\left\langle {\begin{array}{l}
 6 \\
 1 \\
 \end{array}} \right\rangle } \right\rangle =(2)\left\langle {\left\langle
{\begin{array}{l}
 5 \\
 1 \\
 \end{array}} \right\rangle } \right\rangle +\left( {10} \right)\left\langle
{\left\langle {\begin{array}{l}
 5 \\
 0 \\
 \end{array}} \right\rangle } \right\rangle=114
\]
\[
\left\langle {\left\langle {\begin{array}{l}
 7 \\
 1 \\
 \end{array}} \right\rangle } \right\rangle =(2)\left\langle {\left\langle
{\begin{array}{l}
 6 \\
 1 \\
 \end{array}} \right\rangle } \right\rangle +\left( {12} \right)\left\langle
{\left\langle {\begin{array}{l}
 6 \\
 0 \\
 \end{array}} \right\rangle } \right\rangle =240
\]
\[
\left\langle {\left\langle {\begin{array}{l}
 8 \\
 1 \\
 \end{array}} \right\rangle } \right\rangle =(2)\left\langle {\left\langle
{\begin{array}{l}
 7 \\
 1 \\
 \end{array}} \right\rangle } \right\rangle +\left( {14} \right)\left\langle
{\left\langle {\begin{array}{l}
 7 \\
 0 \\
 \end{array}} \right\rangle } \right\rangle =\textbf{494}
\]
\vspace{1mm}

At the end, to be complete, we give besides the general formula of
the second-order Eulerian numbers, which puts into evidence the
connection with the binomial coefficients and with Stirling's
numbers, a representation of the second order Eulerian triangle
(Fig. 6):

\begin{equation}
\label{eq10}
\left\langle {\left\langle {\begin{array}{l}
 n \\
 m \\
 \end{array}} \right\rangle } \right\rangle
=\sum\limits_{k=0}^m {\left( {\begin{array}{l}
 2n+1 \\
 \;\;\;k \\
 \end{array}} \right)} \left\{ {\begin{array}{l}
 n+m+1-k \\
 \;\;m+1-k \\
 \end{array}} \right\}\cdot \left( {-1} \right)^k\quad
for\;n\succ m\ge 0
\end{equation}
where for the binomial coefficients it is valid: $\left( {\begin{array}{l}
 n \\
 k \\
 \end{array}} \right)=\frac{n\;!}{k\;!\;\;\left( {n-k} \right)\;!}$

\noindent and for the Stirling's numbers : $\left\{
{\begin{array}{l}
 n \\
 k \\
 \end{array}} \right\}=k\cdot \left\{ {\begin{array}{l}
 n-\mbox{1} \\
 \;\;k \\
 \end{array}} \right\}+\left\{ {\begin{array}{l}
 n-\mbox{1} \\
 k-\mbox{1} \\
 \end{array}} \right\}\;$.

\centerline{\epsfig{figure=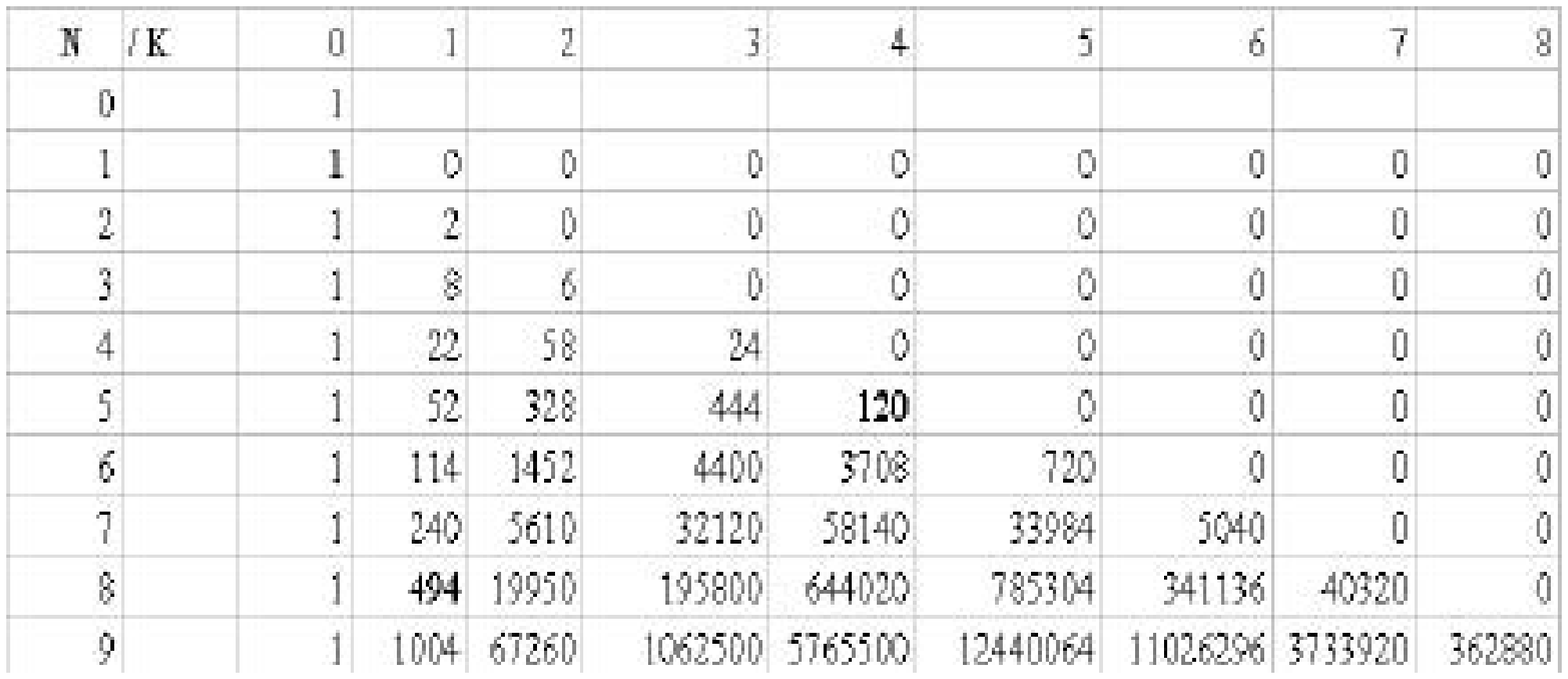,width=11.5cm,height=5cm}}
\begin{center}
{\bf Fig. 6: Second Order Eulerian Triangle }
\end{center}
\baselineskip=0.20in

If we consider, for fundamental of the special function <^s>, the
number $\left\langle {\left\langle {\begin{array}{l}
 {\kern 1pt}\,8 \\
 {\kern 1pt}{\kern 1pt}\,1 \\
 \end{array}} \right\rangle } \right\rangle $ from (\ref{eq10}) we have:
$$\left\langle {\left\langle {\begin{array}{l}
 {\kern 1pt}\,8 \\
 {\kern 1pt}{\kern 1pt}\,1 \\
 \end{array}} \right\rangle } \right\rangle =\left( {\begin{array}{l}
 17 \\
 0 \\
 \end{array}} \right)\;\left\{ {\begin{array}{l}
 10 \\
 2 \\
 \end{array}} \right\}-\left( {\begin{array}{l}
 17 \\
 1 \\
 \end{array}} \right)\;\left\{ {\begin{array}{l}
 9 \\
 1 \\
 \end{array}} \right\}$$
\noindent but keeping in mind that for Stirling's numbers are
valuable the following identities for $n\succ 0$:
\[ \left\{ {\begin{array}{l}
 n \\
 2 \\
 \end{array}} \right\}=2^{n-1}  -1 \;\;\mbox{and}\;\;
\left\{ {\begin{array}{l}
 n \\
 1 \\
 \end{array}} \right\}=1\]
\noindent we finally have:
$$\left\langle {\left\langle {\begin{array}{l}
 8 \\
 1 \\
 \end{array}} \right\rangle } \right\rangle =1\cdot \left(
{2^9-1} \right)-17\cdot 1=511-17=\textbf{494}$$
\medskip

\vspace{2mm}
\begin{center}
5. FROM DISCRETE TO COMPLEX FIELD
\end{center}

\vspace{2mm}

In the previous paragraphs we have discussed on arguments by
infinitesimal, asymptotic, numerical and combinatorial analysis to
characterize the special function <^s> .

It is known that in these fields the completely monotonic functions play a
fundamental role.

We recall that a function $f:I\to \Re $ is said to be completely
monotonic (c.m.) on a real interval $I$, if $f$ has derivatives of all
orders on $I$ which alternate successively in sign, that is:
\[ \;\left( {-1} \right)^n\cdot f^{\left( n \right)}\left( x \right)\ge 0\quad
\forall x\in I \;\;\mbox{and}\;\; \forall \;n\ge 0 \;\;\mbox{with}\;\;
n=0,1,2,3,...\]

In the recent past, various authors\footnote{pag. 445-460, [1].} showed that
numerous functions, which are defined in terms of \textit{gamma},
\textit{polygamma} and other special functions, as the hypergeometrical
ones, are completely monotonic and used this fact to derive many
interesting new inequalities.

We shall confine ourselves only to prove, in very simple way, that
the special function <^s> too, even if it is for definition a
piecewise continuous function, in the real field, it possesses the
same property (c.m.).

\vspace{3mm}

\textbf{Lemma.} \textit{If }$f\left( x \right)\;$\textit{and
}$g\left( x \right)\;$\textit{ are c.m., then }$a\cdot f\left( x
\right)+b\cdot g\left( x \right)\;$\textit{, where a and b}
\textit{ are non-negative constants and }$f\left( x \right)\cdot
g\left( x \right)\;$\textit{are also c.m..}

\vspace{1mm}

The proof of the first thesis is obvious, the second one is then
easily seen from the Leibniz formula :

$$\frac{d^n}{dx^n}\cdot \left[ {f\left( x \right)\cdot g\left( x
\right)} \right]\;=\;\sum\limits_{k=0}^n {\left( {\begin{array}{l}
 n \\
 k \\
 \end{array}} \right)} \;f^{\left( k \right)}\left( x \right)\cdot g^{\left(
{n-k} \right)}\left( x \right)$$

\vspace{2mm}

\textbf{Theorem 1}. \textit{ The special function } <^s> $ \left[
{k,\Omega \;(I_\ell )} \right]\;$\textit{ is completely monotonic in
each }$I_\ell \;$.

\textit{Proof}: The base function $f\left( k
\right)=1+\frac{1}{3k-\Omega \;(I_\ell )}$ is c.m. in each $I_\ell
\;$; in fact the n-th derivative of this function is:

$$\frac{d^n}{dk^n}f\left( k \right)=\frac{\left( {-1}
\right)^n\cdot 3^n\cdot n\,!}{\mathop {\left[ {3k-\Omega \;(I_\ell
)} \right]\,\,}\nolimits^{n+1}} \quad \Rightarrow \quad \left( {-1}
\right)^n\cdot f^{\left( n \right)}\left( k \right)\ge 0\quad
$$

\vspace{1mm}

Thus by the Lemma, in case $k$ is an integer, it's obvious that also
the special function <^s> is c.m. in each $I_\ell \;$; therefore it
remains to prove Theorem 1 for the case when $k$ is a real and
positive number\footnote{ We must also remember the (\ref{eq2})
without giving up the characterization, determined for each interval
$I_\ell \;$. }.

To obtain this it is necessary to the use the following obvious
Theorem 2, which is a consequence of the Lemma, for composed
functions.

\vspace{2mm}

\textbf{Theorem 2: } \textit{Let }$y=f\left( x \right)$
\textit{c.m. and let the power series }$\varphi \left( y
\right)=\sum\limits_{j=0}^\infty {a_j y^j} $\textit{converge for
all }$y$\textit{ in the range of the function }$y=f\left( x
\right). $\textit{ If }$a_j \ge 0$ \textit{for
all}$\;j=0,1,2,3,...$\textit{ then }$\varphi \left[ {f\left( x
\right)} \right]$ \textit{is c.m.}.

\vspace{3mm}

\textbf{Corollary:} \textit{If }$f\left( x \right)$ \textit{is
c.m., then }$e^{f\left( x \right)}$ \textit{is c.m.}.

\vspace{3mm}

In particular, as the special function <^s>  is equal:

$$\mbox{<^s>}\left[ {k,\Omega \;(I_\ell )} \right]\; \;\;=e^{\left( {2k+1}
\right)\cdot \log \left( {1+\frac{1}{3k-\Omega \;(I_\ell )}}
\right)}
$$ it is c.m. in each $I_\ell \;$$. $

\vspace{2mm}

An interesting exposition of the main results on completely monotonic
functions is given in Widder's\footnote{ [7]} work.

That being stated, with some limitations due to the nature of the
special function <^s> in the real field, we can describe a method
for estimating the same function in the complex field with an
important\textit{ improper integral}.

The <^s>  function results to be a piecewise continuous function
because of the presence of the step function $\Omega \;(I_\ell )$,
that is discontinuous, and actually its complete monotonicity has
been proved for all the closed intervals $I_\ell \;$.

Passing to the interval $I\;\left[ {0,\infty \left. \right)}
\right.$ the continuity is not guaranteed and therefore the
application of the following Hausdorff-Bernstein-Widder's\footnote{
Chapter IV, pag. 160-161, [7].} Theorem, must be done carefully, or
it can be limited to characterize the behaviour of the special
function <^s> at the origin and infinity.

\vspace{3mm}

\textbf{Theorem 3}: \textit{A necessary and sufficient condition
for the function }$f\left( s \right)\;$\textit{ in order to be }
\textit{completely monotonic in the interval } $I\;\left[
{0,\infty \left. \right)} \right.$\textit{ is that:}
\begin{equation}
\label{eq11}
f\left( s \right)\;=L_s \left[ {F\left( t \right)}
\right]=\int\limits_0^\infty {e^{-st}dF\left( t \right)}
\end{equation}
\textit{where }$F\left( t \right)$\textit{ is non-decreasing and the
  integral converges in the interval }
$I\;\left[ {0,\infty \left. \right)} \right.$.

\vspace{3mm}

The (\ref{eq11}) represents the transformation of
Laplace-Stieltjes of a locally and absolutely continuous function,
with real values, in the interval $I\;\left[ {0,\infty \left.
\right)} \right.$.

Now, keeping in mind the relation existing in such case between
Laplace-Stieltjes transform $L_s $ and \textit{ordinary}\footnote{ Chapter
  IV, [2].} Laplace transform $L$ ( observing that we can suppose
$F\left( 0 \right)=0)$:
\begin{equation}
\label{eq12}
L_s \left[ {F\left( t \right)} \right]=s\cdot \int\limits_0^\infty
{e^{-st}F\left( t \right)\cdot dt} =s\cdot L\left[ {F\left( t \right)}
\right]
\end{equation}
we can determine the expression of the function $F\left( t \right)$ by the
inverse Laplace transform.

From (\ref{eq11}) and (\ref{eq12}) we get
$\Phi \left( s \right)= \frac{f\left( s\right)}{s}= L\left[
{F\left( t \right)} \right]$ and successively we need
to face the calculus of integral of the type\footnote{ pag. 168, [2].}:

\vspace{4mm}

\qquad \qquad \qquad \qquad $F\left( t \right)\mbox{
}=\mbox{L}^{\mbox{-1}}\left[ {\Phi \left( s \right)}
\right]=v.p.\,\frac{1}{2\pi \,i}\cdot \int\limits_{x_0 -i\cdot
\infty }^{x_0 +i\cdot \infty } {e^{ts}\Phi \left( s \right)\cdot
ds} \qquad \qquad \qquad  (14.1)$

\vspace{4mm}

\noindent with $t\succ 0\,,\,\Phi \left( s \right)$ holomorphic
function in the half-plane $\Re e\left( s \right)\succ 0$ and $x_0 $
arbitrary real positive number.
\newpage
In general the calculus of Bromwich's integral (14.1), that is
defined as a \textit{ Cauchy principal value}, we can only do it
numerically, applying quadrature formulae, but in our case [
$f\left( s \right)\;=\mbox{<^s> } \left( s \right)\;$], being the
$\Phi \left( s \right)$ function \textbf{piecewise analytic}, we
could prove by the direct calculus of the considered integral (see
the following paragraph), that $F\left( t \right)=2\cdot \Theta
\left[ t \right]$ and $\Theta \left[ t \right]$ represents the unit
step function or Heaviside's function (\,$\Theta \left[ {\;t<0\;}
\right]=0\;;\;\Theta \left[ {\;t>0\;} \right]=1$ ) and in
(\ref{eq11}) we'll have $dF\left( t \right)=2\cdot d\,\Theta \left[
t \right]=2\cdot \delta \left( t \right)dt$, with $\delta \left( t
\right)$ that is the distribution of Dirac\footnote{ pag. 30-35,
[2].}.

\vspace{6mm}
\setnashbf
\begin{center}
6. THE BEHAVIOUR OF THE SPECIAL FUNCTION <^s> AT THE ORIGIN AND
INFINITY\footnote{ pag. 491-492, [4].}
\end{center}

\setnash

\vspace{2mm}

In the previous paragraph we have stated the following
approximation ``$\approx $'', in terms of Laplace transform:

\[ \Phi \left( s \right)= \quad \mbox{<^s> } \left( s \right)\;\raise0.7ex\hbox{}
\!\mathord{\left/ {\vphantom {
s}}\right.\kern-\nulldelimiterspace}\!\lower0.7ex\hbox{s}\approx L\left[
{F\left( t \right)} \right]=L\left[ {2\cdot \Theta \left( t \right)}
\right]=\int\limits_0^\infty e^{-st}2\cdot \Theta \left( t
\right)\,dt=\frac{2}{s}\]

The approximation is essentially origined by neglecting the point of
discontinuities of the first kind of the special function <^s> ,
between an interval $I_\ell \;$and the following $I_{\ell +1} \;$as
far as the interval $I\;\left[ {0,\infty \left. \right)} \right.$.

Said that to calculate the Bromwich's integral (14.1), we consider
the path of integration rightly deformed, as we can see in the Fig.
7.

\medskip
\centerline{\epsfig{figure=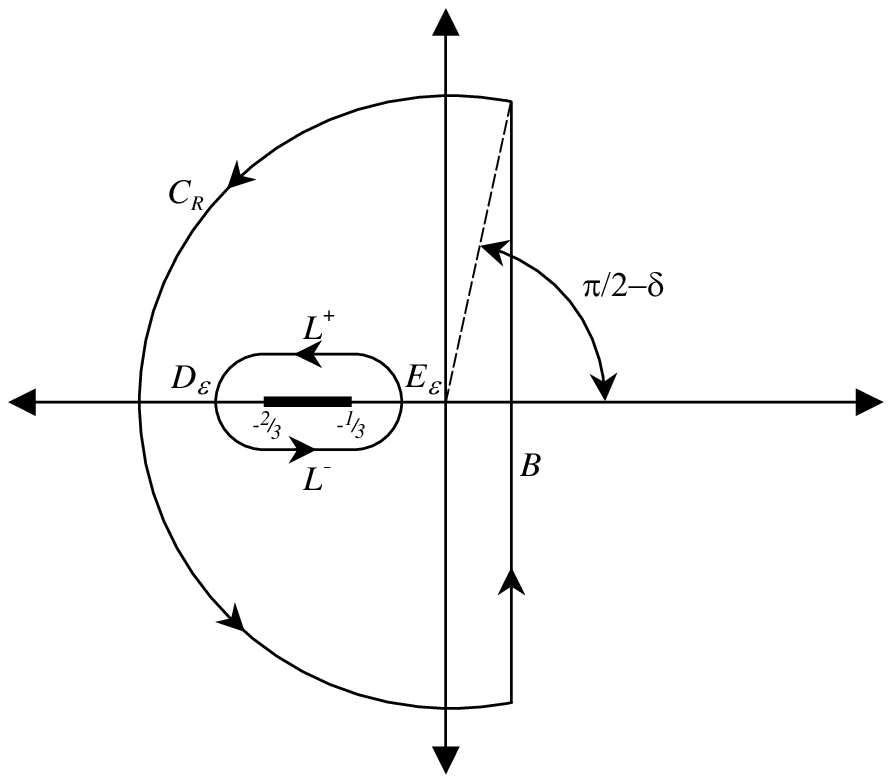,width=8cm,height=7cm}}
\medskip
\baselineskip=0.20in
\begin{center}
{\bf Fig. 7: The path of integration}
\end{center}

$\Phi \left( s \right)$ is a multiple-valued analytic function; in
the interest of obtaining a single value mapping we consider the
principal branch of the power in $\mbox{<^s> } \left( s \right)\;$,
that allows us to treat a holomorphic branch of the same $\Phi
\left( s \right)$ function.

In its field of single-values, the function $\Phi \left( s \right)$
possesses a simple pole at $s=0\;$and two branch points at
$s=-\frac{1}{3}\;$ and $s=-\frac{2}{3}\;$; let's observe, in fact,
that in the interval $\left[ {-1,0} \right]$ the special function
<^s> is:
$\mbox{<^s> } \left( s \right)\;=\left( {1+\frac{1}{3s+1}}
\right)^{2s+1}\;\;=e^{\left( {2s+1} \right)\cdot \log \left(
{1+\frac{1}{3s+1}} \right)}$

A cut, joining the two branch points, would prevent $s$ to circulate
around them, and the special function $\mbox{<^s> } \left( s
\right)\;$can be treated as a \textbf{piecewise holomorphic
function}.

In conclusion, defined with $C$ the boundary (shown in Fig. 7), with $C_R $
the circular arc of radius $R$, with $B$ the vertical line with
$\Re \left( s\right)=x_0 $, with $L$ the boundary of the cut branch, formed by
$D_\varepsilon $ and $E_\varepsilon $, that is the semi-circles of radius
$\varepsilon $ capping the ends of the branch cut and by $L^+$ and $L^-$,
lines above and below, we'll have:
\begin{equation}
\label{eq13}
\frac{1}{2\pi \,i}\cdot \int\limits_{x_0 -i\cdot \infty }^{x_0 +i\cdot
\infty } {e^{ts}\Phi \left( s \right)\cdot ds} =
\frac{1}{2\pi \,i}\left\{ \, \right.\int\limits_C {...} -\int\limits_{C_R }
{...} -\int\limits_{L^+} {...} -\int\limits_{L^-} {...}
-\int\limits_{D_\varepsilon } {...} -\int\limits_{E_\varepsilon } {...}
\left.  \right\}
\end{equation}

The first integral in the second member of (\ref{eq13}), by the
residue theorem, is $2\cdot \Theta {\kern 1pt}\left( t \right)$, as
the first order pole, at the origin, gives: $\lim_{s\to 0}
\mbox{<^s> } \left( s \right)\;=2$.

The second integral with $s=R\cdot e^{i\vartheta }$ is:
\begin{equation}
\label{eq14}
\int\limits_{C_R } {e^{ts}\Phi \left( s \right)\cdot ds}
=\int\limits_{\pi \mathord{\left/ {\vphantom {\pi 2}} \right.
\kern-\nulldelimiterspace} 2-\delta }^{\pi \mathord{\left/ {\vphantom {\pi
2}} \right. \kern-\nulldelimiterspace} 2} {...d\vartheta } +\int\limits_{\pi
\mathord{\left/ {\vphantom {\pi 2}} \right. \kern-\nulldelimiterspace}
2}^{{3\pi } \mathord{\left/ {\vphantom {{3\pi } 2}} \right.
\kern-\nulldelimiterspace} 2} {...d\vartheta } +\int\limits_{{3\pi }
\mathord{\left/ {\vphantom {{3\pi } 2}} \right. \kern-\nulldelimiterspace}
2}^{{3\pi } \mathord{\left/ {\vphantom {{3\pi } 2}} \right.
\kern-\nulldelimiterspace} 2+\delta } {...d\vartheta }
\end{equation}
and therefore it vanishes: in fact the first and third integral in the
second member of (\ref{eq14}) vanish as $R \to \infty $ by the maximum
modulus bound and the second integral vanished by Jordan's Lemma.

The third and fourth integral in the second member of (\ref{eq13}) cancel each
other along the paths $L^+$ and $L^-$: in fact their values, calculated
along their opposite paths, eliminate each other.

In the end, the last two integrals of (\ref{eq13}), by the maximum modulus
bound, vanish as $\varepsilon  \to 0$.

Now, from the asymptotic behaviour of the function $F\left( t
\right)$, so determined, we can therefore deduce asymptotic
properties of the correspondent Laplace transform, that is to use
the following Abelian theorems (initial and final value theorem):

\vspace{4mm}

\textbf{Theorem:} \textit{Let }$F$\textit{ be a transformable
function and let's suppose that the} $\lim_{t\to \infty } F\left(
t \right)$ \textit{ exists, then the } $\lim_{s\to 0}\;\left[
\right.s\cdot \Phi \left( s \right)\left. \right]$ \textit{
exists, too (let's suppose, for convenience }$s\in \Re $\textit{)
  and is:}
\begin{equation}
\label{eq15}
\lim _{s\to 0} \;\left[ \right.s\cdot \Phi \left( s \right)\left. \right]=
\lim _{t\to \infty } F\left( t \right)
\end{equation}
\indent \textit{If the }$\lim _{t\to 0^+} F\left( t
\right)$\textit{ exists, then the }$\lim_{s\to \infty } s\cdot
\Phi \left(s \right)$\textit{ exists, too and is (with }$s\in \Re
):$
\begin{equation}
\label{eq16}
\lim _{s\to \infty } \;\left[ \right.s\cdot \Phi \left( s \right)\left.
\right]=
\lim _{t\to 0^+} F\left( t \right)
\end{equation}

In our case we have seen that $F\left( t \right)=2\cdot \Theta {\kern
1pt}\left( t \right)$ , with $\Theta {\kern 1pt}\left( t \right)$, that is
Heaviside's function, and therefore we have that:
\[
\lim _{t\to \infty } F\left( t \right) =\lim _{t\to \infty }
2\cdot \Theta {\kern 1pt}\left( t \right)=2
\]
and also:
$$\lim _{t\to 0^+} F\left( t \right)=\lim _{t\to 0^+} 2\cdot \Theta
{\kern 1pt}\left( t \right)=2$$

The existence and the calculus of such limits, in (\ref{eq15}) and
(\ref{eq16}), give the following results (between them a further
confirm of the \textit{fundamental} theorem):

\[ \lim _{s\to 0}  \quad \mbox{<^s> } \left( s \right)\;=2\mbox{ ; }
\lim _{s\to \infty } \mbox{<^s> } \left( s \right)\;=2\]

\newpage

\begin{center}
\textbf{APPENDIX}
\end{center}

\textbf{\textit{DERIVE Version 6 : }}\textbf{the vector function }:
S(k,$\Omega$) := $\lim _{\;t\to k} \;\left( {1+\frac{1}{3t-\Omega}}
\right)^{2t+1}$

\vspace{4mm}

\begin{center}
{\#}1: V(k, $\Omega$ ) := VECTOR([t, o, S(t,$\Omega$ )], t, k, k +
10)
\end{center}

\vspace{4mm}

\begin{center}
\textbf{\textit{a display result}}\textbf{\textit{:}}
\end{center}

\vspace{8mm}

V($ $ k,\quad $\Omega$ )  \qquad\qquad\qquad\qquad\qquad    VECTOR

$\vert \quad 09\quad 1\quad
2.0484148121729077984789748528464903812082646338028 \quad \vert $

$\vert \quad 10\quad 1\quad
2.0379259208387064562838079920238964441117176933744 \quad \vert $

$\vert \quad 11\quad 1\quad
2.0294161672236677191636908626945714916029532064944 \quad \vert $

$\vert \quad 12\quad 1\quad
2.0223737073469397533461445484297949184415534016560 \quad \vert $

$\vert \quad 13\quad 1\quad
2.0164491799135882361365114303236301480590762568587 \quad \vert $

$\vert \quad 14\quad 1\quad
2.0113959747189663594458436566806886371655211482846 \quad \vert $

$\vert \quad 15\quad 1\quad
2.0070350457364044054984268130457357906298862812402 \quad \vert $

$\vert \quad 16\quad 1\quad
2.0032332566108411453651981972386971733090123101528 \quad \vert $

$\vert \quad 17\quad 1\quad
1.9998895526624551656968593976078763119133058370586 \quad \vert $

$\vert \quad 18\quad 1\quad
1.9969258468076576081148471529242828805292201418873 \quad \vert $

$\vert \quad 19\quad 1\quad
1.9942808454379732420411582337304540085201659981375 \quad \vert $

\newpage
\baselineskip=0.20in

\textit{Kragujevac J. Math.} 29 (2006) 26--35.

\vspace{20mm}

\baselineskip=0.30in

\setnashbf
\begin{center}
\textbf{THE SPECIAL FUNCTION } <^s> , \textbf{II.}
\end{center}
\setnash
\begin{center}
\vspace{10mm}

{\large \bf  Andrea Ossicini}

\vspace{10mm}

\baselineskip=0.20in

{\it Via delle Azzorre 352-D2, 00121 Roma, Italy.\\
{\rm (e-mail: a.ossicini@finsiel.it)}

\vspace{2mm}

(submitted september 9, 2005) }

\end{center}

\vspace{2mm}

\baselineskip=0.17in \noindent {\small {\bf Abstract.} We describe a
method for estimating the special function <^s> , in the complex cut
plane $A = {\mathbf{C}}\backslash \;\left( { - \infty ,0} \right]$,
with a Stieltjes transform, which implies that the function <^s> is
\textit{logarithmically completely monotonic}. To be complete, we
find a nearly exact integral representation. At the end, we also
establish that  $1 \mathord{\left/ {\vphantom {1 }} \right.
\kern-\nulldelimiterspace} \mbox{<^s>} \left( x \right)$ is a
complete Bernstein function and we give the representation formula
which is analogous to the L\'evy-Khinchin formula.

\vspace{3mm}

2000\textit{ Mathematics Subject Classification: primary }33B15;
\textit{secondary }26A48; 33E20.

\vspace{2mm}

Keywords: special functions; completely monotonic functions;
integral transforms; Bernstein functions.

\baselineskip=0.20in

\vspace{3mm}

\begin{center}
1. INTRODUCTION
\end{center}

\vspace{2mm}

In [8] the author introduces a new special function, named with the
Arabian letter\footnote{ The letter <^s> (shin) is the thirteenth
letter of the Arabian alphabet.} <^s>, and proves that this is
logarithmically convex and completely monotonic for all the closed
real intervals $I_\ell \;$with $\ell = 1,2,3,...$ .

\vspace{1mm}

The explicit formula of the special function <^s> , in the discrete
field, is:

\begin{equation}
\label{eq19} \mbox{<^s>}\left[ {k,\Omega \;(I_\ell )} \right]\;
\quad = \;\left( {1 + \frac{1}{3k - \Omega \;(I_\ell )}\mbox{ }}
\right)^{2k + 1}\;
\end{equation}

\vspace{1mm}

\noindent where is always valid \textit{the following
boundary}\footnote{ The boundary can include the sign ``='' if the
integer variable $k$ goes towards zero or the infinity.}:

\vspace{2mm}

$\mbox{<^s>} \left[ {k,\Omega \;(I_\ell -1)} \right] \prec 2\prec
 \mbox{<^s>} \left[ {k,\Omega \;(I_\ell )} \right]\;\;$ with
$\Omega \;(I_\ell -1)=\Omega \;(I_\ell )-1$ ; $\forall \;I_\ell
\;,\;k,\ell \in \mathbf{N}$

\newpage

The \textit{auxiliary }integer function $\Omega(I_\ell )$, that
really represents a growing ``step function'', is defined, for the
intervals of 8 or 9 following values of $k$, in the following way:

\begin{itemize}
\item $\Omega \;(I_1 )=0$ for $k$=1,{\ldots},8
\item $\Omega \;(I_2 )=1$ for $k$=9,{\ldots},16 ; $\Omega \;(I_3 )$=2 for $k$=17,{\ldots},25 ; $\Omega \;(I_4 )$=3 for $k$=26,{\ldots},34
\item $\Omega \;(I_5 )=4$ for $k$=35,{\ldots},43 ; $\Omega \;(I_6 )$ =5 for $k$=44,{\ldots},51 ; $\Omega \;(I_7 )$=6 for $k$=52,{\ldots},60
\item $\Omega \;(I_8 )=7$ for $k$=61,{\ldots},69; $\Omega \;(I_9 )$ =8 for $k$=70,{\ldots},78 ; $\Omega \;(I_{10} )$=9 for $k$=79,{\ldots},86
\item $\Omega \;(I_{11} )=10$ for $k$=87,{\ldots},95 ; $\Omega \;(I_{12} )$=11 for $k$=96,{\ldots},104. ; etc.
\end{itemize}

Successively we give (Fig. 1) the graphs, related to the families of
<^s> functions, that are $ \mbox{<^s>} \left[ {k,\Omega \;(I_\ell )}
\right]\;$ and $ \mbox{<^s>} \left[ {k,\Omega \;(I_\ell - 1)}
\right]\;$, or better, to the set of the arcs belonging to them, and
to the auxiliary function $\Omega \;(I_\ell )$ .

\vspace{6mm}
\medskip
\centerline{\epsfig{figure=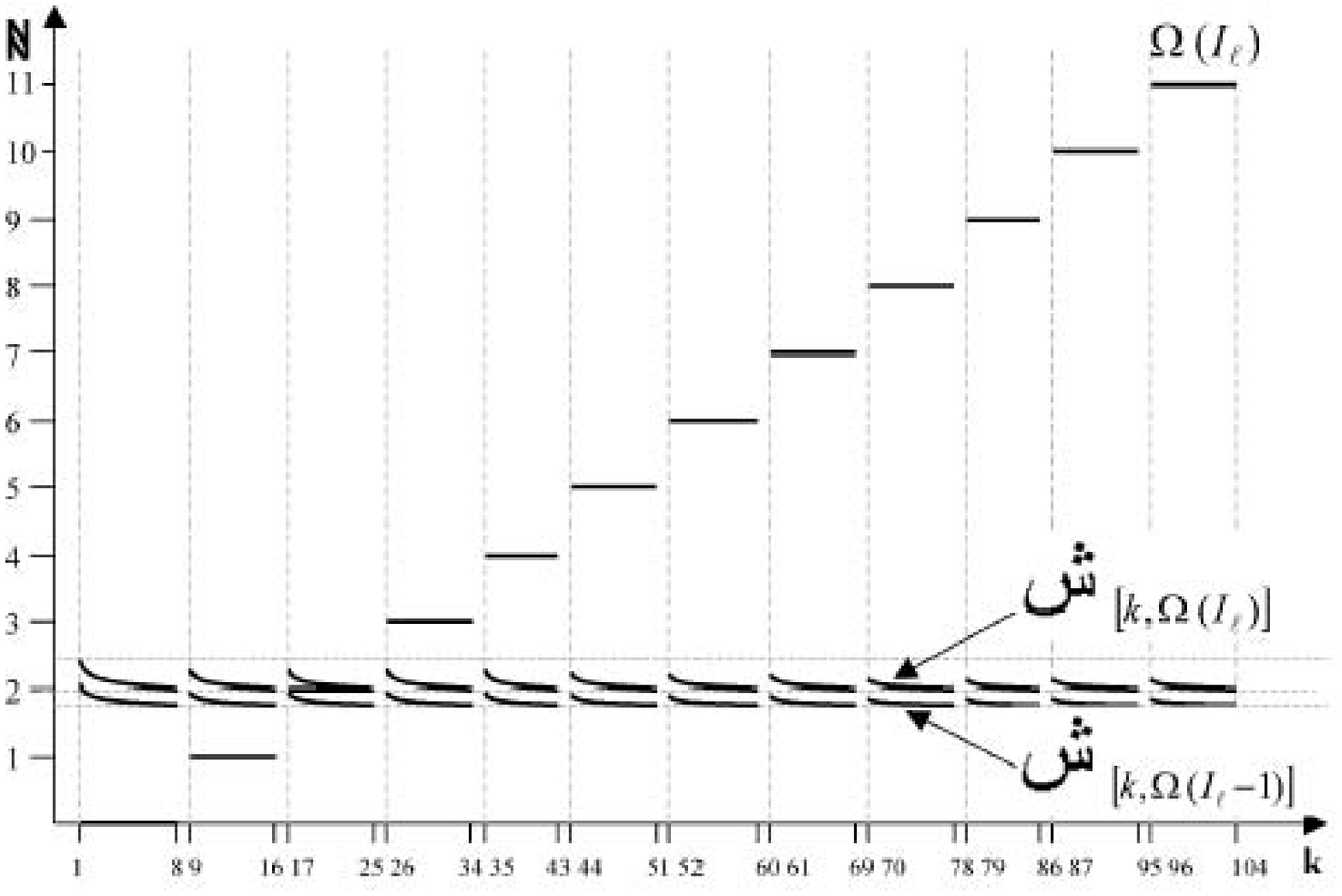,width=13cm,height=8.5cm}}
\medskip
\baselineskip=0.20in
\begin{center}
{\bf Fig. 1.}
\end{center}

Let's extend the dependence of the integer step function $\Omega
\;(I_l )$ to the real field and let's use the following definition:

\vspace{3mm}

\[ \Omega\left( x\right)=\mbox{min}\left\{
{k\in}\right.\mathbf{N}\mbox{: S}_{k+1}\left( x\right)\ge
2\left.\right\}\mbox{; }x\in\mathbf{R^{+}}\mbox{ and where S}_k
\left( x\right)=\left({1+\frac{1}{3x-k+1}} \right)^{2x+1}\]

\vskip 0.1truecm

By simple algebraic passages we have, by a more appropriate
notation, due to Iverson, that:
\begin{equation}
\label{eq20} \Omega\;\left( x \right)=\left\lceil
{3x-\frac{1}{2^{\frac{1}{2x+1}}-1}} \right\rceil
\end{equation}
where $\left\lceil {\,x} \right\rceil $ means the smallest integer,
greater than $x$ or equal to it.

\vskip 0.1truecm

\vspace{3mm}

Therefore the special function <^s> $ $ possesses the following
explicit formula in the real field:

\[ \mbox{<^s>}\left( x \right)=\left(
{1+\frac{1}{3x-\Omega \left( x \right)}} \right)^{2x+1}=\left(
{1+\frac{1}{3x-\left\lceil {3x-\frac{1}{2^{\frac{1}{2x+1}}-1}}
\right\rceil }} \right)^{2x+1}\]\

Extending the field of definition of the variable $k$ to the real
positive numbers, it's possible to notice that such function, being
represented by the union of continuous arcs (all above the straight
line of height 2, see Fig.1) is actually assimilable to a piecewise
continuous function.

\vspace{3mm}

That being stated, an important subclass of completely monotonic
functions consists of the Stieltjes transforms defined as the class
of functions $f$: $\left( {0,\infty } \right) \to {\rm R}$ of the
form:

\begin{equation}
\label{eq21} f\left( x \right)\; = a + \int\limits_0^\infty
{\frac{d\mu \;\left( t \right)}{x + t}}
\end{equation}

\noindent where $a \ge 0$ and $\mu \left( t \right)$ is a
nonnegative measure on $\;\left[ {0,\infty \left. \right)} \right.$
with $\int\limits_0^\infty {\frac{d\mu \;\left( t \right)}{1 + t}}
\le \infty$, see [2].

\vspace{3mm}

In the Addenda and Problems in ([1], p.127), it is stated that if a
function $f$ is holomorphic in the cut plane $A =
{\mathbf{C}}\backslash \;\left( { - \infty ,0} \right]$ and
satisfies the following conditions :

\vspace{3mm}

 (i) $\Im f\left( z \right) \le 0$ for $\Im \left( z
\right) \succ 0$

\vspace{2mm}

(ii) $f(x) \ge 0$ for $x \succ 0$

\vspace{3mm}

\noindent then $f$ is a Stieltjes transform.

\newpage

\begin{center}
2. THE REPRESENTATION AS A STIELTJES TRANSFORM
\end{center}

\vspace{4mm}

In ([8], \S 6) the author characterizes the holomorphy
(\textbf{piecewise analytic}) of the special function $\mbox{<^s>}
\left( z \right)\;$ in the cut plane $B = {\mathbf{C}}\backslash
\,\left[ { - \raise0.7ex\hbox{2} \!\mathord{\left/ {\vphantom {2
3}}\right.\kern-\nulldelimiterspace}\!\lower0.7ex\hbox{3}, -
\raise0.7ex\hbox{1} \!\mathord{\left/ {\vphantom {1
3}}\right.\kern-\nulldelimiterspace}\!\lower0.7ex\hbox{3}} \right]$
and proves a remarkable result that implies :

\begin{equation}
\label{eq22} \lim  _{\left| z \right| \to \infty } \mbox{<^s>} (z) =
2 \quad \quad (\mbox{z} \in B)
\end{equation}

\vspace{3mm}

To prove that the harmonic function $\; \Im (\mbox{<^s>})\;$
satisfies $\;\Im \mbox{<^s>} \left( z \right) \le 0\;$ for $\;\Im
\left( z \right) \succ 0\;$, we use that maximum principle for
subharmonic functions, that can be found in ([4], p. 20), and show
that lim sup of $\;\Im (\mbox{<^s>})\;$ at all boundary points
including infinity is less than or equal to 0.

\vspace{3mm}

From (\ref{eq22}) we conclude that this is true at infinity.

\vspace{3mm}

How, for definition $\mbox{<^s>} \left( x \right) \succ 0\;$ for $x
\succ 0$; these last statements imply the result (\ref{eq21}).

\vspace{3mm}

The constant $a$ in (\ref{eq21}) is given by :

\[
a = \lim _{x \to \infty } \mbox{<^s>} \left( x \right)\;
\]

\noindent and therefore for the \textit{fundamental }theorem of the
special function  <^s> $ $ we have, see ([8], \S 2):

\[
a = \lim _{x \to \infty } \mbox{<^s>} \left( x \right)  = 2
\]

In (\ref{eq21}) $\mu \left( t \right)$ is the limit in the vague
topology of measures

\[
d\mu \left( t \right) = \lim _{y \to 0^ + } \;\; - \frac{1}{\pi }
\Im f\left( { - t + iy} \right)dt
\]

For $z \in B = {\mathbf{C}}\backslash \left[ { - \raise0.7ex\hbox{2}
\!\mathord{\left/ {\vphantom {2
3}}\right.\kern-\nulldelimiterspace}\!\lower0.7ex\hbox{3}, -
\raise0.7ex\hbox{1} \!\mathord{\left/ {\vphantom {1
3}}\right.\kern-\nulldelimiterspace}\!\lower0.7ex\hbox{3}} \right]$
we have in the close interval $\left[ { - 1,0} \right]$ :

\vspace{3mm}

\begin{equation}
\label{eq23} \mbox{<^s>} \left( z \right)\; =  \left( {1 +
\frac{1}{3z + 1}} \right)^{2z + 1}\; = \mbox{ }\exp \left( {\left(
{2z + 1} \right) \cdot Log\left( {1 + \frac{1}{3z + 1}} \right)}
\right)
\end{equation}

\vspace{3mm}

\noindent where \textit{Log }denotes the principal branch of the
logarithm.

\newpage

Let $t \in {\mathbf{R}}$ and $z \in {\mathbf{C}}$ with $\Im \left( z
\right) \succ 0$.

\vspace{2mm}

If $z$ tends to $t$ , then for (\ref{eq19}) and (\ref{eq23}),
results ( with $\ell \in {\mathbf{Z}}$ ) \footnote{ ${\mathbf{Z}}$
denotes the relative integer set.}:

\[
\mbox{<^s>} \left( z \right)\; = \left\{ {\begin{array}{l}
 \;\left( {1 + \frac{1}{3t - \Omega \;(I_\ell )}} \right)^{2t + 1}\quad
\quad if \;t \succ 0 \\
 \;\;\;\; \\
 \quad \;\quad \;2\quad \quad \quad  \;if\;\; t = 0\; \\
 \; \\
 \;\exp \left[ {\left( {2t + 1} \right)\, \cdot \log \left( {\left|
{\frac{3t + 2}{3t + 1}} \right|} \right)\, - \,k \cdot i\,\pi \left(
{2t + 1} \right)} \right] \; with\;\;k = 1\;\;\; if \; - \frac{2}{3}
\prec \,t\, \prec - \frac{1}{3}
\\
 \\
 \;\;\; and\;\; with\;\;k = 0\;\;\; if\;- 1 \prec \,t\, \prec - \frac{2}{3}
\;\;\;
or\;\;\;if \; - \frac{1}{3} \prec t \prec 0\ \\
 \\
 \quad \;\quad \;2\quad \quad \quad \;if\;\; t = - 1\; \\
 \;\;\; \\
 \;\; \\
 \;\left( {1 + \frac{1}{3t - \Omega \;(I_\ell )}} \right)^{2t + 1}\quad
\;if\;t \prec - 1 \\
 \end{array}} \right. \\
 \\
\]
\vspace{3mm}

In particular then we obtain , if y tends to $0^ + , $
for$\footnote{ In the discontinuity points 1/3 and 2/3 we
respectively compute the limits of the real variable t on the left
and on the right (see Fig. 2).}$ $t \in {\mathbf{R}}$:

\[
 - \frac{1}{\pi }\Im f\left( { - t + iy} \right)\quad \to
\left\{ {\begin{array}{l}
 \quad 0\quad if\;\;t \le \raise0.7ex\hbox{$1$} \!\mathord{\left/ {\vphantom
{1 3}}\right.\kern-\nulldelimiterspace}\!\lower0.7ex\hbox{$3$}\quad
or\quad t \ge \raise0.7ex\hbox{$2$} \!\mathord{\left/ {\vphantom {2
3}}\right.\kern-\nulldelimiterspace}\!\lower0.7ex\hbox{$3$} \\
 \\
 \\
 \quad \frac{1}{2\pi }\frac{\left( {\left( {3t - 1} \right)^2}
\right)\,^t}{\left( {\left( {3t - 2} \right)^2} \right)\,^t} \cdot
\sin \;\left( {2\pi t} \right) \cdot \left\{ {\left| {\frac{3t -
2}{3t - 1}} \right| - \frac{3t - 2}{3t - 1}} \right\}\quad
\;\;if\;\;\raise0.7ex\hbox{$1$} \!\mathord{\left/ {\vphantom {1
3}}\right.\kern-\nulldelimiterspace}\!\lower0.7ex\hbox{$3$} \prec t
\prec \raise0.7ex\hbox{$2$} \!\mathord{\left/ {\vphantom {2
3}}\right.\kern-\nulldelimiterspace}\!\lower0.7ex\hbox{$3$} \\
 \end{array}} \right.\;
\]

\noindent and using the identity (the Euler reflection formula):

\[
\Gamma \left( \alpha \right) \cdot \Gamma \left( {1 - \alpha }
\right) = \frac{\pi }{\sin \left( {\alpha \pi } \right)}
\]

\noindent we are now in a position to determine the following nearly
exact integral representation ( Stieltjes transform):

\begin{equation}
\label{eq24} \mbox{<^s>} \left( x \right)\;
 \approx  2 + \frac{1}{2} \cdot \;\int\limits_{\raise0.7ex\hbox{$1$}
\!\mathord{\left/ {\vphantom {1
3}}\right.\kern-\nulldelimiterspace}\!\lower0.7ex\hbox{$3$}}^{\raise0.7ex\hbox{$2$}
\!\mathord{\left/ {\vphantom {2
3}}\right.\kern-\nulldelimiterspace}\!\lower0.7ex\hbox{$3$}}
{\,\left\{ {\frac{1}{\Gamma \left( {2 \cdot t} \right) \cdot \Gamma
\left( {1 - 2 \cdot t} \right)}\frac{\left( {\left( {3t - 1}
\right)^2} \right)\,^t}{\left( {\left( {3t - 2} \right)^2}
\right)^{\,t}}\left[ {\left| {\frac{3t - 2}{3t - 1}} \right| -
\frac{3t - 2}{3t - 1}} \right]} \right\}\,\frac{dt}{\left( {x + t}
\right)}}
\end{equation}

The approximation is essentially origined by neglecting the point of
discontinuities of the first kind of the special function <^s> $ $,
in the real field, between an interval $I_\ell \;$ and the following
$I_{\ell +1} \;$ as far as the interval $I\;\left[ {0,\infty \left.
\right)} \right.$.
\newpage

Successively we give (Fig. 2: $x \to t$) the graphs (\textit{red
color}) related to the function\footnote{ $\mathop {\mathop \mu
\limits^\bullet } \left( t \right) = \frac{d\mu \left( t
\right)}{dt}$}  $\mathop {\mathop \mu \limits^\bullet } \left( t
\right)$ in the real interval $I\left[ {-1,1} \right]$: for $t =
\raise0.7ex\hbox{$1$} \!\mathord{\left/ {\vphantom {1
2}}\right.\kern-\nulldelimiterspace}\!\lower0.7ex\hbox{$2$}
\Rightarrow$ $\mathop {\mathop \mu \limits^\bullet } \left( t
\right) = 0$ and this point is a \textit{flex point} with oblique
tangent.

\vspace{4mm}

\medskip
\centerline{\epsfig{figure=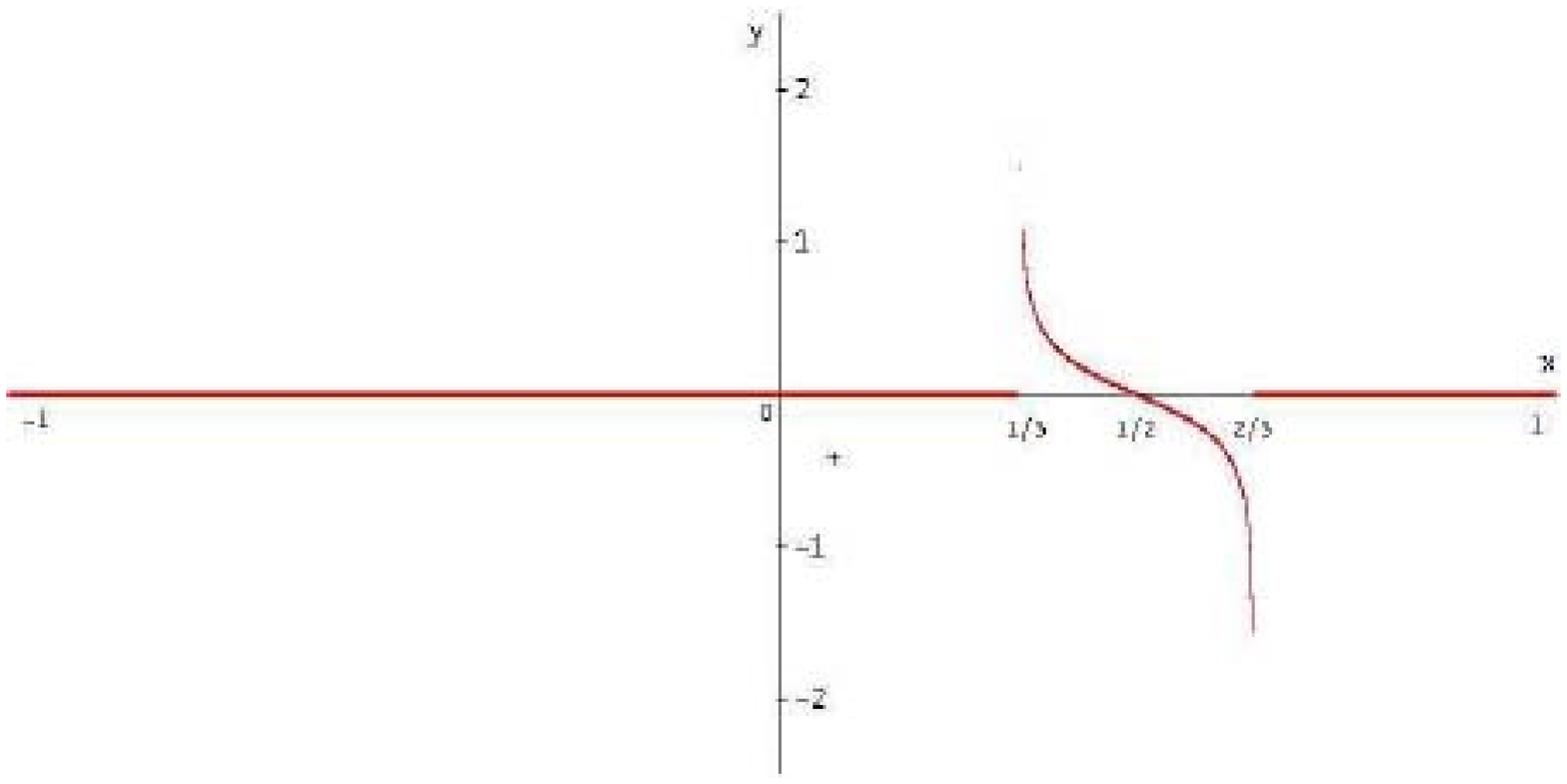,width=16cm,height=8cm}}
\medskip
\baselineskip=0.20in
\begin{center}
{\bf Fig. 2: the graph of} $\mathop {\mathop \mu \limits^\bullet }
\left( t \right)$
\end{center}

\vspace{4mm}

That being stated, we denote the set of completely monotonic
functions with $ \mathcal{C}$.

\vspace{2mm}

Now, we also recall that a function $f$: $\left] {0,\infty } \right[
\to \left] {0,\infty } \right[$ is said  to be logarithmically
completely monotonic [5], if it is $C^\infty $ and

\[
\;\left( { - 1} \right)^k \cdot \left[ {\log f\left( x \right)}
\right]^{\left( k \right)} \ge 0\quad \quad for\;k = 1,2,3,...
\]

To simplify we denote the class of logarithmically completely
monotonic functions by $\mathcal{L}$ and the set of Stieltjes
transforms by $\mathcal{S}$ .

\vspace{2mm}

In order to prove that the special function $ \mbox{<^s>} \left( x
\right)$ is logarithmically completely monotonic, we need the
following \textbf{lemma :}

\begin{center}
$\mathcal{S}$\textbf{\textit{ $\backslash $ }}$\left\{ 0 \right\}
\subset $ $\mathcal{L}$
\end{center}

\newpage

This lemma is a consequence of the following result, established by
Horn [6], that allows also to characterize the class of
logarithmically completely functions as the infinitely divisible
completely monotonic functions:

\vspace{6mm}

\textbf{Theorem 1: } \textit{For a function} $f$: $\left] {0,\infty
} \right[ \to \left] {0,\infty } \right[$  \textit{the following are
equivalent:}

\vspace{4mm}
\begin{center}
(i) $f \in $ $\mathcal{L}$ ; \quad (ii) $f^\alpha \in $
$\mathcal{C}$ for all $\alpha \succ 0$ and $\alpha \in {\mathbf{R}}$
; \quad (iii) $\sqrt[n]{f} \in $ $\mathcal{C}$ for all $n =
1,2,3...$
\end{center}

In fact, let $f \in $ $\mathcal{S}$  ( $\mathcal{S} \subset
\mathcal{C}$ ) and non-zero and let $\alpha \succ 0$, by Theorem 1
it is immediate to prove that $f^\alpha \in \mathcal{C}$.

\vspace{4mm}

Now, writing $\alpha = n + a$ with $n = 0,1,2,...$ and $0 \le a
\prec 1$ we have $f^\alpha = f^n \cdot f^a$, and using the stability
of $\mathcal{C}$\textbf{\textit{ }}under multiplication and that
$f^a \in $\textbf{\textit{ }}$S$ $ \Rightarrow  \quad
S$\textbf{\textit{ $\backslash $ }}$\left\{ 0 \right\} \subset $
$\mathcal{L}$.

\vspace{4mm}

In conclusion for (\ref{eq24}) also the special function
$\mbox{<^s>} \left( x \right) \in  \mathcal{L}$.

\vspace{7mm}

\begin{center}
3. THE CLASS OF BERNSTEIN FUNCTIONS
\end{center}

\vspace{4mm}

There is an important relation between the set $\mathcal{S} $ of
Stieltjes transforms and the class $\mathcal{B} $ of Bernstein
functions.

\vspace{3mm}

We recall that a function $f$ : $\left( {0,\infty } \right) \to
\left[ {0,\infty } \right)$ is called a Bernstein function, if $f$
 has derivatives of all orders and $f'$ is completely monotonic.

\vspace{3mm}

Now, if $f$ is non-zero Stieltjes transform, then  $1
\mathord{\left/ {\vphantom {1 f}} \right. \kern-\nulldelimiterspace}
f$  is a Bernstein function  ([3], Prop. 1.3).

\vspace{3mm}

The special function $\mbox{<^s>} \left( x \right) \in $
$S$\textbf{\textit{ $\backslash $ }}$\left\{ 0 \right\} $ and this
fact implies that $1 \mathord{\left/ {\vphantom {1 }} \right.
\kern-\nulldelimiterspace} \mbox{<^s>} \left( x \right)$ is a
Bernstein function.

\vspace{3mm}

In addition, using the identity :

\[
\raise0.7ex\hbox{$1$} \!\mathord{\left/ {\vphantom {1
{\mbox{<^s>}\left( x
\right)}}}\right.\kern-\nulldelimiterspace}\!\lower0.7ex\hbox{${\mbox{<^s>}
\left( x \right)}$} = \raise0.7ex\hbox{$x$} \!\mathord{\left/
{\vphantom {x {x \cdot \mbox{<^s>}\left( x
\right)}}}\right.\kern-\nulldelimiterspace}\!\lower0.7ex\hbox{${x
\cdot \mbox{<^s>}\left( x \right)}$}
\]

\noindent and remembering the following definition:

\vspace{2.5mm}

A Bernstein function $\phi $ is called a special Bernstein function
if the function $\raise0.7ex\hbox{$\lambda $} \!\mathord{\left/
{\vphantom {\lambda {\phi \left( \lambda
\right)}}}\right.\kern-\nulldelimiterspace}\!\lower0.7ex\hbox{${\phi
\left( \lambda \right)}$}$ is also a Bernstein function.

\vspace{2.5mm}

\noindent we can conclude that  $ \;\; x\cdot\mbox{<^s>}\left( x
\right)\; \; $ is a special Bernstein function.

\vspace{2.5mm}

The family of special Bernstein functions is very large, and it
contains in particular the family of complete Bernstein functions
(also known as \textit{operator-monotone functions}, see [7], for
instance).

\vspace{2.5mm}

Recall that a function $\phi :\left( {0,\infty } \right) \to \Re $
is called a complete Bernstein function if there exits a Bernstein
function $\eta $ such that :

\[
\phi (\lambda ) = \lambda ^2L\left[ {\eta \left( \lambda \right)}
\right] \quad , \quad \lambda \succ 0
\]

\noindent where $L$ stands for the Laplace transform.

\vspace{2.5mm}

Now, using the main results about the special function $\mbox{<^s>}
\left( x \right)$ ([8], \S 5 and \S 6) it is immediate to establish
that $x \cdot \mbox{<^s>}\left( x \right)$ is a complete Bernstein
function.

\vspace{2.5mm}

Note also that a function $f\left( x \right)$ is called a complete
Bernstein function if, and only if,

\begin{equation}
\label{eq25}\ f\left( x \right) = a + bx + \int\limits_{0 + }^\infty
{\frac{x}{t + x}} \rho \left( {dt} \right) \
\end{equation}

\noindent where $a,b \ge 0$ and $ \rho $ is a Radon measure on
$\left({0,\infty}\right)$ such that $\int\limits_{0 + }^\infty
{\left( {1/(1 + t)}\right)} \rho \left( {dt} \right) < \infty$ .

\vspace{2.5mm}

From this one, we may deduce that the function $\; x \to {f\left( x
\right)} \mathord{\left/ {\vphantom {{f\left( x \right)} x}} \right.
\kern-\nulldelimiterspace} x\; $ is a Stieltjes transform [2].

\vspace{2.5mm}

This result was actually already obtained with the representation
(\ref{eq24}) of Stieltjes of the special function <^s>.

\vspace{2.5mm}

At the end, recall that the following conditions are equivalent:

\vspace{4mm}

(i) $\phi $ is a complete Bernstein function;

\vspace{2.5mm}

(ii) $\raise0.7ex\hbox{$\lambda $} \!\mathord{\left/ {\vphantom
{\lambda {\phi \left( \lambda
\right)}}}\right.\kern-\nulldelimiterspace}\!\lower0.7ex\hbox{${\phi
\left( \lambda \right)}$}$ is a complete Bernstein function.

\vspace{3mm}

This result implies also that the function $1 \mathord{\left/
{\vphantom {1 }} \right. \kern-\nulldelimiterspace} \mbox{<^s>}
\left( x \right)$ is a complete Bernstein function and, remembering
the standard form (\ref{eq25}) and that the functions $1
\mathord{\left/ {\vphantom {1 }} \right. \kern-\nulldelimiterspace}
\left[ x \cdot\mbox{<^s>} \left( x \right)\right]$ and $1
\mathord{\left/ {\vphantom {1 }} \right. \kern-\nulldelimiterspace}
\mbox{<^s>} \left( \frac{1}{x} \right)$ are Stieltjes transforms, it
is easy and immediate to estabilish that the constant $a$ ({\it
killing rate}) and $b$ ({\it drift coefficient}) are given by:

\[
a = \lim _{x \to \infty } 1 \mathord{\left/ {\vphantom {1 }} \right.
\kern-\nulldelimiterspace} \mbox{<^s>} \left( \frac{1}{x} \right)  =
1 \mathord{\left/ {\vphantom {1 }} \right.
\kern-\nulldelimiterspace} \mbox{<^s>} \left( 0 \right)\; =
 \frac{1}{2}\;
\]

\[
b = \lim _{x \to \infty } 1 \mathord{\left/ {\vphantom {1 }} \right.
\kern-\nulldelimiterspace} \left[ x \cdot\mbox{<^s>} \left( x
\right)\right] = 0 \;
\]

\vspace{2.5mm}

\noindent and the following representation formula which is
analogous to the L\'evy-Khinchin formula:

\vspace{2.5mm}

\[ 1 \mathord{\left/
{\vphantom {1 }} \right. \kern-\nulldelimiterspace} \mbox{<^s>}
\left( x \right)\;
 \approx \frac{1}{2} + \frac{1}{2} \cdot \int\limits_{0 + }^\infty
{\,\left\{ {\frac{1}{\Gamma \left( {2 \cdot t} \right) \cdot \Gamma
\left( {1 - 2 \cdot t} \right)}\frac{\left( {\left( {3t - 2}
\right)^2} \right)\,^t}{ \left( {\left( {3t - 1} \right)^2}
\right)^{\,t}}\left[ {\left| {\frac{3t - 1}{3t - 2}} \right| -
\frac{3t - 1}{3t - 2}} \right]} \right\}\,\frac{x\cdot t }{\left( {t
+ x} \right)}} \; {dt}
\]

\vspace{3mm}

For the interplay between complete Bernstein functions and Stieltjes
transforms we refer also to [9].

\vspace{3mm}

Finally, with an Euler-Venn diagram, we give the most important
analytic properties of the special function <^s> .

\vspace{3mm}

\medskip
\centerline{\epsfig{figure=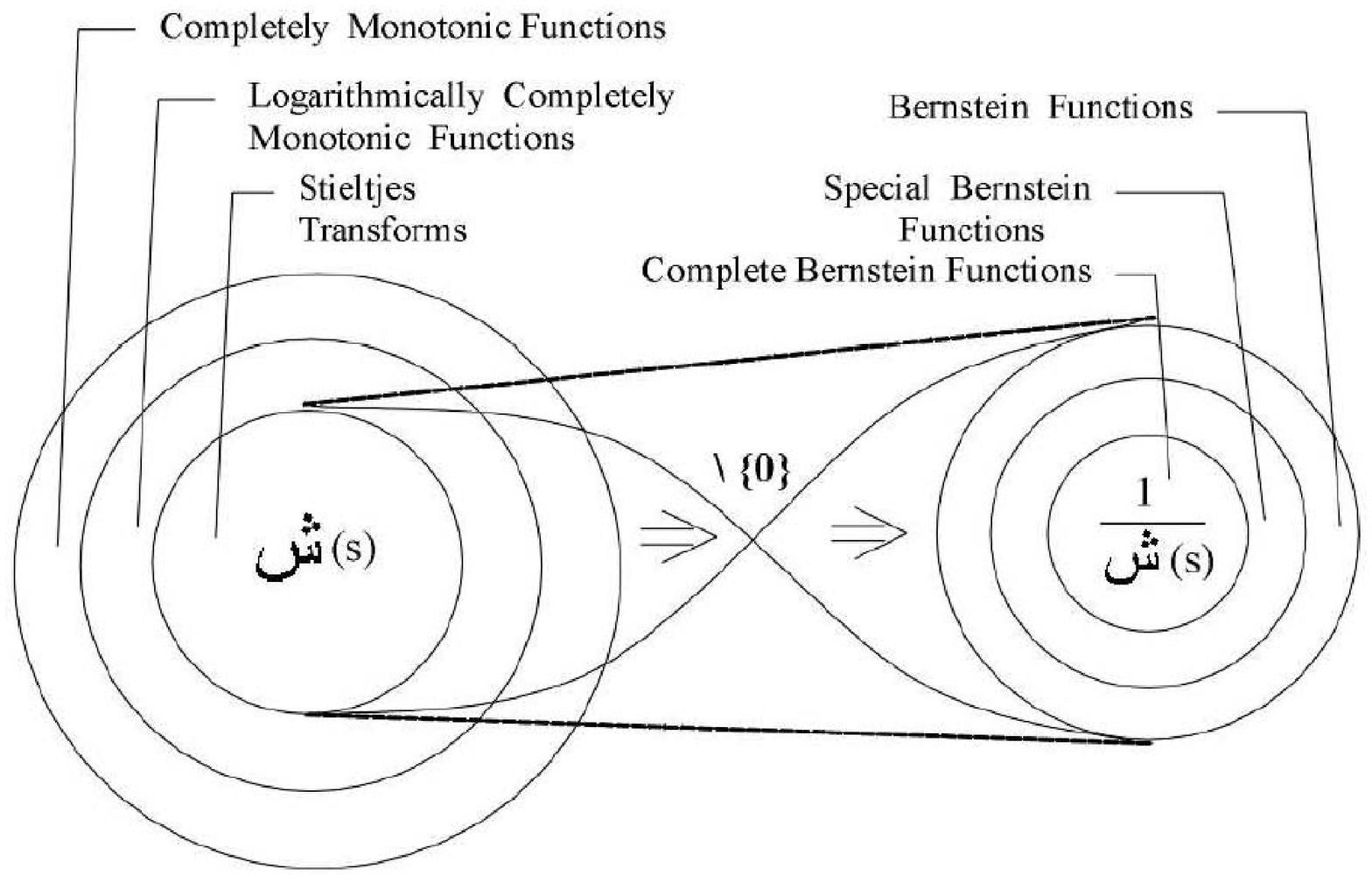,width=15cm,height=10cm}}
\medskip
\baselineskip=0.20in \baselineskip=0.20in
\begin{center}
{\bf Completely Monotonic Functions $\sim$ vs  $\sim$ Bernstein
Functions}
\end{center}

\newpage

 \baselineskip=0.20in
\textit{Kragujevac J. Math.} 31 (2007) 36--45.

\vspace{20mm}

\setarab

\baselineskip=0.30in

\setnashbf
\begin{center}
\textbf{THE SPECIAL FUNCTION } <^s> , \textbf{III.}
\end{center}
\setnash
\begin{center}
\vspace{10mm}

{\large \bf  Andrea Ossicini}

\vspace{10mm}

\baselineskip=0.20in

{\it Via delle Azzorre 352-D2, 00121 Roma, Italy.\\
{\rm (e-mail: a.ossicini@finsiel.it)}

\vspace{4mm}

(submitted February 03, 2006) }

\end{center}

\vspace{4mm}

\baselineskip=0.17in \noindent {\small {\bf Abstract.} We prove the
following exact symbolic formula of the special function  <^s>, $ $
in the entire s-complex plane with the negative real axis (including
the origin) removed, with a double Laplace transform:

\vspace{6mm}

$\;\;\;\;\;\;\;\;\;\;\;\;\;\;\;\;\mbox{<^s>}\;(s) = L\left\{ \right.
2 \cdot \delta \left( t \right) + L \left\{ \right. \frac{1}{2\pi
i}\cdot $[$\;\mbox{<^s>}\;(t \cdot e^{\; - i\pi }) - \mbox{<^s>}\;
(t \cdot e^{\;i\pi })$]$\left. \right\} \left. \right\} $

\vspace{6mm}

\noindent where $\delta \left( t \right)$ stands for the
distribution of Dirac and $e$ represents the Euler's number.

\baselineskip=0.20in

\vspace{4mm}

\begin{center}
1. THE EXACT SYMBOLIC FORMULA
\end{center}

\vspace{2mm}

In [4] and [5] the author introduces and characterizes the discrete
and special function  <^s> $ $ in the complex field.

\vspace{1mm}

We recall that the explicit formula of the special function <^s> ,
in the discrete field, is:

\[ \mbox{<^s>}\left[ {k,\Omega \;(I_\ell )} \right]\; \quad = \;\left(
{1 + \frac{1}{3k - \Omega \;(I_\ell )}\mbox{ }} \right)^{2k + 1}\]

\noindent where it is always valid \textit{the following
boundary}\footnote{ The boundary can include the sign ``='' if the
integer variable $k$ goes towards zero or the infinity.}:

\vspace{3mm}

$\mbox{<^s>} \left[ {k,\Omega \;(I_\ell -1)} \right] \prec 2\prec
 \mbox{<^s>} \left[ {k,\Omega \;(I_\ell )} \right]\;\;$ with
$\Omega \;(I_\ell -1)=\Omega \;(I_\ell )-1$ ; $\forall \;I_\ell
\;,\;k,\ell \in \mathbf{N}$

\newpage

The \textit{auxiliary }integer function $\Omega(I_\ell )$, that
really represents a growing ``step function'', is defined, for the
intervals of 8 or 9 following values of $k$, in the following way:

\begin{itemize}
\item $\Omega \;(I_1 )=0$ for $k$=1,{\ldots},8
\item $\Omega \;(I_2 )=1$ for $k$=9,{\ldots},16 ; $\Omega \;(I_3 )$=2 for $k$=17,{\ldots},25 ; $\Omega \;(I_4 )$=3 for $k$=26,{\ldots},34
\item $\Omega \;(I_5 )=4$ for $k$=35,{\ldots},43 ; $\Omega \;(I_6 )$ =5 for $k$=44,{\ldots},51 ; $\Omega \;(I_7 )$=6 for $k$=52,{\ldots},60
\item $\Omega \;(I_8 )=7$ for $k$=61,{\ldots},69; $\Omega \;(I_9 )$ =8 for $k$=70,{\ldots},78 ; $\Omega \;(I_{10} )$=9 for $k$=79,{\ldots},86
\item $\Omega \;(I_{11} )=10$ for $k$=87,{\ldots},95 ; $\Omega \;(I_{12} )$=11 for $k$=96,{\ldots},104. ; etc.
\end{itemize}

\vspace{4mm}

In particular the author proves in ([4], \S \S 5-6) that the special
function <^s> $ $ in the real field is completely monotonic for all
the closed intervals $I_\ell \;$ and therefore the following first
approximation ``$ \approx $'' in the complex field:

\begin{equation}
\label{eq26} \mbox{<^s>} \left( s \right) \approx L_s \left[
{F\left( t \right)} \right] = L_s \left[ {2 \cdot \Theta \left( t
\right)} \right]
 =
\int\limits_0^\infty e^{ - st}d\left[ {2 \cdot \,\Theta \left( t
\right)} \right]\;= 2
\end{equation}

\noindent where $L_s $ denotes the Laplace-Stieltjes transform and
$\Theta \left[ t \right]$ represents the unit step function or
Heaviside's function (\,$\Theta \left[ {\;t<0\;}
\right]=0\;;\;\Theta \left[ {\;t>0\;} \right]=1$ )

\vspace{3mm}

Now, we observe that the Laplace-Stieltjes transform is closely to
related other integral transforms, including the Fourier transform
and the Laplace transform.

\vspace{3mm}

In particulary if g has derivative g', then the Laplace-Stieltjes
transform of g is the Laplace transform of g'.

\vspace{3mm}

Consequently, considering the derivative of Heaviside step function,
from (\ref{eq26}) we  have :

\begin{equation}
\label{eq27} \mbox{<^s>}\left( s \right) \approx L \left[ {2 \cdot
\delta \left( t \right)} \right] = \int\limits_0^\infty e^{ -
st}2\delta \left( t \right)\mbox{dt} = 2
\end{equation}

\noindent where $L$ denotes the ordinary Laplace transform and
$\delta \left( t \right)$ stands for distribution of Dirac.

\vspace{3mm}

In ([5], \S 2) the author find practically the following nearly
exact integral representation:

\begin{equation}
\label{eq28} \mbox{<^s>}\left( x \right)
 \approx 2 + \frac{1}{2}\int\limits_0^\infty {\,\left\{ {\frac{1}{\Gamma
\left( {2 \cdot t} \right) \cdot \Gamma \left( {1 - 2 \cdot t}
\right)}\frac{\left( {\left( {3t - 1} \right)^2} \right)\,^t}{\left(
{\left( {3t - 2} \right)^2} \right)^{\,t}}\left[ {\left| {\frac{3t -
2}{3t - 1}} \right| - \frac{3t - 2}{3t - 1}} \right]}
\right\}\,\frac{dt}{\left( {x + t} \right)}}
\end{equation}

\vspace{3mm}

\noindent where $\Gamma$ denotes the \textit{Eulerian gamma
function} (or the \textit{Eulerian integral of second kind}).

\vspace{3mm}

The approximation is essentially origined by neglecting the point of
discontinuities of the first kind of the special function <^s> $ $,
in the real field, between an interval $I_\ell \;$ and the following
$I_{\ell +1} \;$ as far as the interval $I\;\left[ {0,\infty \left.
\right)} \right.$.

\vspace{3mm}

We recall that extending the field of definition of the integer
variable $k$ of the discrete and special function  <^s> $ $ to the
real positive number,  it's possible to notice that such function is
actually assimilable to a piecewise continuous function.

\vspace{3mm}

On the contrary of equation (\ref{eq27}), in the equation
(\ref{eq28}) there is an Stieltjes transform, that arises naturally
as an iteration of the ordinary Laplace transform.

\vspace{3mm}

In fact, if

\[
f\left( x \right) = \int\limits_0^\infty e^{ - xt}\varphi \,\left( t
\right)\,dt
\]

\noindent where

\[
\varphi \,\left( x \right) = \int\limits_0^\infty e^{ - xt}\psi
\,\left( t \right)\,dt
\]

\noindent then, changing the order of  integrations in the double
integral by appealing to Fubini's Theorem,  we have formally ([1],
p. 127 and [8], p. 335):

\begin{equation}
\label{eq29} f\left( x \right) = L\left\{ {L\,\,\left[ {\psi \,t}
\right]} \right\}\, = \int\limits_0^\infty e^{ -
xu}\,du\;\int\limits_0^\infty e^{ - ut}\psi \left( t \right)\,dt
\end{equation}

\[
 = \int\limits_0^\infty \psi \left( t \right)\,dt\;\int\limits_0^\infty e^{
- u\left( {x + t} \right)}\,du
\]

Hence

\begin{equation}
\label{eq30} f\left( x \right) = \int\limits_0^\infty \frac{\psi
\left( t \right)}{x + t}\,dt
\end{equation}

\vspace{3mm}

This last equation we refer to as Stieltjes transform, or from
another point view as the Stieltjes integral equation.

\vspace{3mm}

However, we shall usually be concerned with the more general case in
which the integral equation (\ref{eq30}) is replaced by a Stieltjes
integral:

\[
f\left( x \right)\; = \int\limits_0^\infty {\frac{d\alpha \;\left( t
\right)}{x + t}}
\]

\noindent and in this form the equation was considered by T.J.
Stieltjes in connection with his work on continued fractions [6].

\newpage

In (\ref{eq28}) the Stieltjes integral, that is analogous to
(\ref{eq30}), converges and then we have $ $ (see [8],Theorem 7b, p.
340):

\begin{equation}
\label{eq31} \lim _{\,\eta \to 0 + } \frac{f\left( { - \xi - i\,\eta
} \right) - f\left( { - \xi + i\,\eta } \right)}{2\pi \;i} =
\frac{\psi \left( {\xi + } \right) + \psi \left( {\xi - }
\right)}{2}
\end{equation}

\vspace{3mm}

\noindent for any positive $\xi $ at which $\psi \left( {\xi + }
\right)$ and $\psi \left( {\xi - } \right)$ exist.

\vspace{3mm}

For, simple computation gives :

\[
\frac{f\left( { - \xi - i\,\eta } \right) - f\left( { - \xi +
i\,\eta } \right)}{2\pi \;i}
 = \frac{1}{\pi }\,\int\limits_0^\infty \frac{\eta \psi \left( t
\right)}{\left( {t - \xi } \right)^2 + \eta ^2}\,dt
\]

The integral is known as Poisson's integral for the half-plane or as
Cauchy's singular integral ([7], p. 30).

\vspace{3mm}

The result (\ref{eq31}) may also be written symbolically as

\begin{equation}
\label{eq32} \frac{f\left( {x \cdot e^{ - i\pi }} \right) - f\left(
{x \cdot e^{i\pi }} \right)}{2\pi \;i} = \psi \left( x \right)
\end{equation}

\vspace{3mm}

In fact it is sufficient to compare the equation 11.8.4 in ([7], p.
318).

\vspace{3mm}

That being stated in ([4], \S 6) the author characterizes the
holomorphy (\textbf{piecewise analytic}) of the special function $
\mbox{<^s>} \left( s \right)\;$ in the cut plane $A =
{\mathbf{C}}\backslash \;\,\left[ { - \raise0.7ex\hbox{$2$}
\!\mathord{\left/ {\vphantom {2
3}}\right.\kern-\nulldelimiterspace}\!\lower0.7ex\hbox{$3$}, -
\raise0.7ex\hbox{$1$} \!\mathord{\left/ {\vphantom {1
3}}\right.\kern-\nulldelimiterspace}\!\lower0.7ex\hbox{$3$}}
\right]$.

\vspace{3mm}

In the close interval $\left[ { - 1,0} \right]$ results:

\begin{equation}
\label{eq33} \mbox{<^s>}\left( s \right)\; = \quad \left( {1 +
\frac{1}{3s + 1}} \right)^{2s + 1}\;\; = \mbox{ } \exp \left(
{\left( {2s + 1} \right) \cdot Log\left( {1 + \frac{1}{3s + 1}}
\right)} \right)
\end{equation}

\vspace{3mm}

\noindent where \textit{Log }denotes the principal branch of the
logarithm in the interest of obtaining a single value mapping.

\vspace{3mm}

In effects the special function <^s> $ $ possesses two branch points
at $s=-\frac{1}{3}\;$ and $s=-\frac{2}{3}\;$ and therefore a cut,
joining the two branch points, would prevent $s$ to circulate around
them, and the special function $\mbox{<^s> } \left( s \right)\;$can
be treated as a holomorphic function.

\vspace{3mm}

We observe that, in our particular case, the results (\ref{eq32})
and (\ref{eq31}) with (\ref{eq33}) can be written and compute in the
following way:

\vspace{3mm}

\[
\frac{\mbox{<^s>} \left( {x \cdot e^{ - i\pi }} \right) -
\mbox{<^s>} \left( {x \cdot e^{i\pi }} \right)}{2\pi i} = \quad \lim
_{y \to 0^ + } \frac{\mbox{<^s>} \left( { - x - iy} \right) -
\mbox{<^s>} \left( { - x + iy} \right)}{2\pi i}
\]

\newpage

$= \lim _{y \to 0^ + }
 \left\{ \frac{\exp \left( {\left( {2\left( {-x-iy} \right) + 1} \right) \cdot
Log\left( {1 + \frac{1}{3\left( {-x-iy} \right) + 1}} \right)}
\right)- \exp \left( {\left( {2\left( {-x+iy} \right) + 1} \right)
\cdot Log\left( {1 + \frac{1}{3\left( {-x+iy} \right) + 1}} \right)}
\right)}{2\pi i} \right\} $

\vspace{3mm}

\begin{equation}
\label{eq34}
 =  \quad \frac{\sin \;\left( {2\pi t} \right)}{2\pi }\cdot \frac{\left( {\left( {3t - 1} \right)^2}
\right)\,^t}{\left( {\left( {3t - 2} \right)^2} \right)\,^t} \cdot
\left\{ {\left| {\frac{3t - 2}{3t - 1}} \right| - \frac{3t - 2}{3t -
1}} \right\}
\end{equation}

\vspace{3mm}

Finally, with (\ref{eq27}),(\ref{eq28}),(\ref{eq29}),(\ref{eq34})
and using the identity (the Euler reflection formula):

\vspace{3mm}

\[
\Gamma \left( \alpha \right) \cdot \Gamma \left( {1 - \alpha }
\right) = \frac{\pi }{\sin \left( {\alpha \pi } \right)}
\]

\vspace{3mm}

\noindent we can prove, for the linearity property of the Laplace
transformation, the following exact symbolic formula of the special
function {<^s>}, in the cut plane $A = {\mathbf{C}}\backslash
\;\left( { - \infty ,0} \right]$:

\vspace{3mm}

\begin{equation}
\label{eq35} \mbox{<^s>}\;(s) = L\left\{ \right. 2 \cdot \delta
\left( t \right) + L \left\{ \right. \frac{1}{2\pi i}\cdot
[\;\mbox{<^s>}\;(t \cdot e^{\; - i\pi }) - \mbox{<^s>}\; (t \cdot
e^{\;i\pi })]\left. \right\} \left. \right\}
\end{equation}

\vspace{3mm}

\noindent where $\delta \left( t \right)$ stands for the
distribution of Dirac and $e$ represents the Euler's number.

\vspace{4mm}

In conclusion we give also the three-dimensional graph  of absolute
value of the special function <^s>  of a complex variable $s = x +
i\,y$ (Fig. 1) from \textit{DERIVE} computer algebra system and in
APPENDIX we explain the choice of the Arabian letter <^s> $ $.

\medskip
\centerline{\epsfig{figure=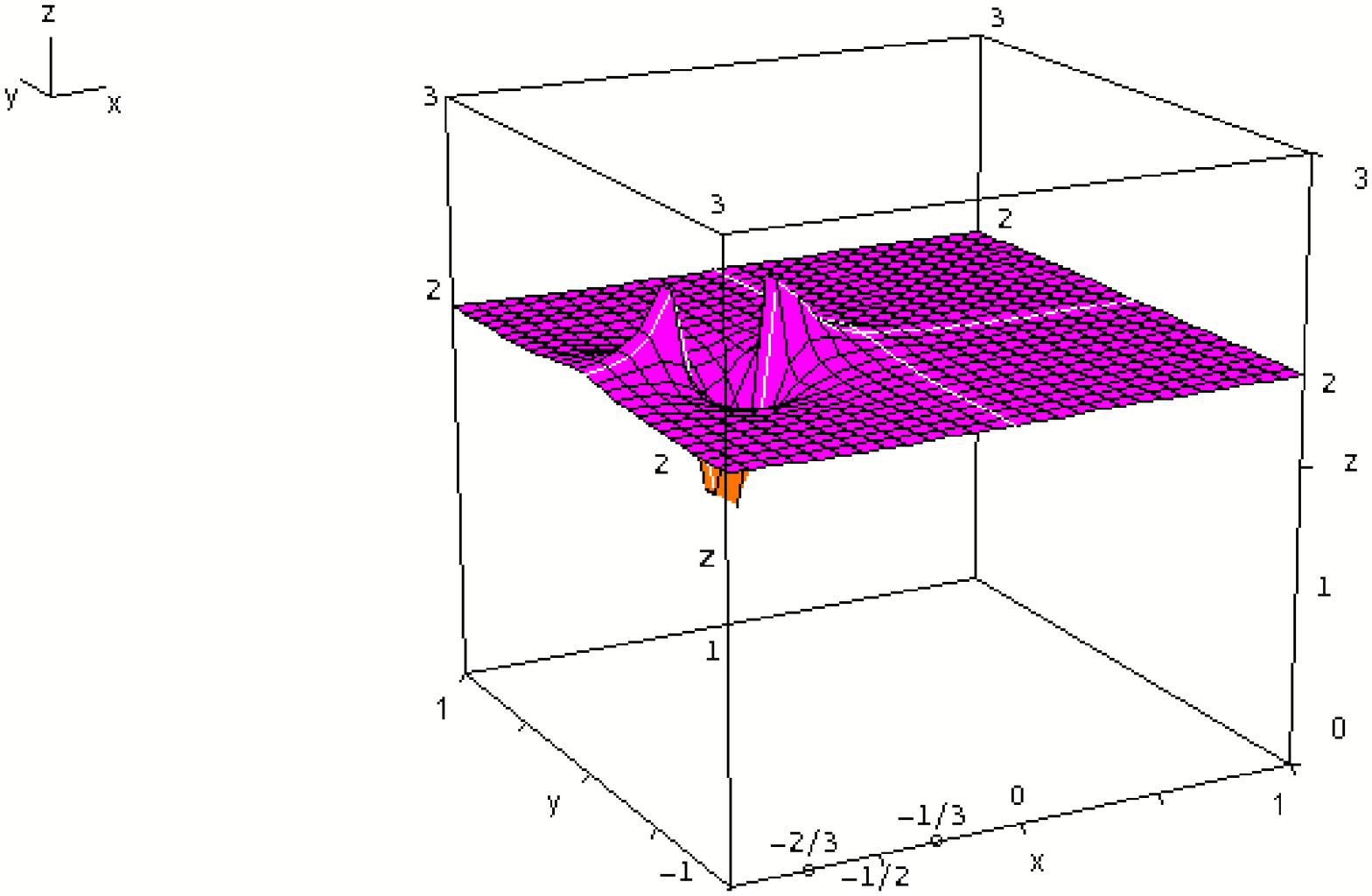,width=17.75cm,height=8cm}}
\medskip
\baselineskip=0.20in
\begin{center}
{\bf Fig. 1}
\end{center}

\newpage

\begin{center}
2. ADDITIONAL ANALYTIC REMARKS
\end{center}

\vspace{4mm}

The beautiful symbolic formula (35) is in practice a consequence of
the fact that Laplace-Stieltjes transform are often written as
ordinary Laplace transform involving the distribution of Dirac,
sometimes referred to as the symbolic impulse function $\delta
\left( t \right)$.

\vspace{3mm}

Besides that the construction-algorithm ([4], \S 1) of the discrete
and special function <^s> $ $ had already evidenced the necessity of
an asymmetric impulsive force in order to maintain the set of the
points belonging to them all above the horizontal straight line of
height 2.

\vspace{3mm}

Now, we observe that the asymmetrical impulse function $\delta _ +
\left( t \right)$ is more suitable for use in connection with the
one-sided Laplace transformation that the symmetrical impulse
function $\delta \left( t \right)$.

\vspace{3mm}

In effects in the formula (35) the presence of the symbolic impulse
function $\delta _ + \left( t \right)$ would be more appropriated ,
also because if one applies the Laplace transformation to the
``definition'' of the impulse function $\delta _ + \left( t
\right)$, one obtains the formal result:

\[
L\left[ {\delta _ + \left( t \right)} \right]\; = \;1
\]

We recall that the asymmetrical impulse function $\delta _ + \left(
t \right)$ is defined ([3], see Sec. 21.9-6) by:

\vspace{2mm}

\[
\begin{array}{l}
 \int_{a + 0}^b {f\left( \xi \right)\;} \delta _ + \left( {\xi - x}
\right)\;d\;\xi \; = \left\{ {_{f\left( {x + 0} \right)\;\quad\;
\mbox{if}\;\;a \le x \prec b\;}^{\mbox{ 0}\quad \quad \quad
\mbox{if}\;\;x \prec
a\;or\;x \ge b} } \right\}\;\quad \left( {a \prec b} \right) \\
 \\
 \end{array}
\]

It is possible to write:

\[
\delta _ + \left( t \right) = \frac{d}{dx}U_ + \left( x \right)
\]

\noindent where $U_ + \left( x \right)$ denotes the asymmetrical
unit-step function :

\[
U_ + \left( x \right)\; = \left\{ {_{\mbox{ 1}\;\quad for\;x \succ
0\;}^{\mbox{ 0}\quad\; for\;x \le 0} } \right\}
\]

\vspace{3mm}

Successively we give the cartesian diagram (Fig. 2) of the unit-step
function $U_ + \left( x \right)$ and of the one related to
approximation of the impulse function $\delta _ + \left( t \right)$.

\medskip
\centerline{\epsfig{figure=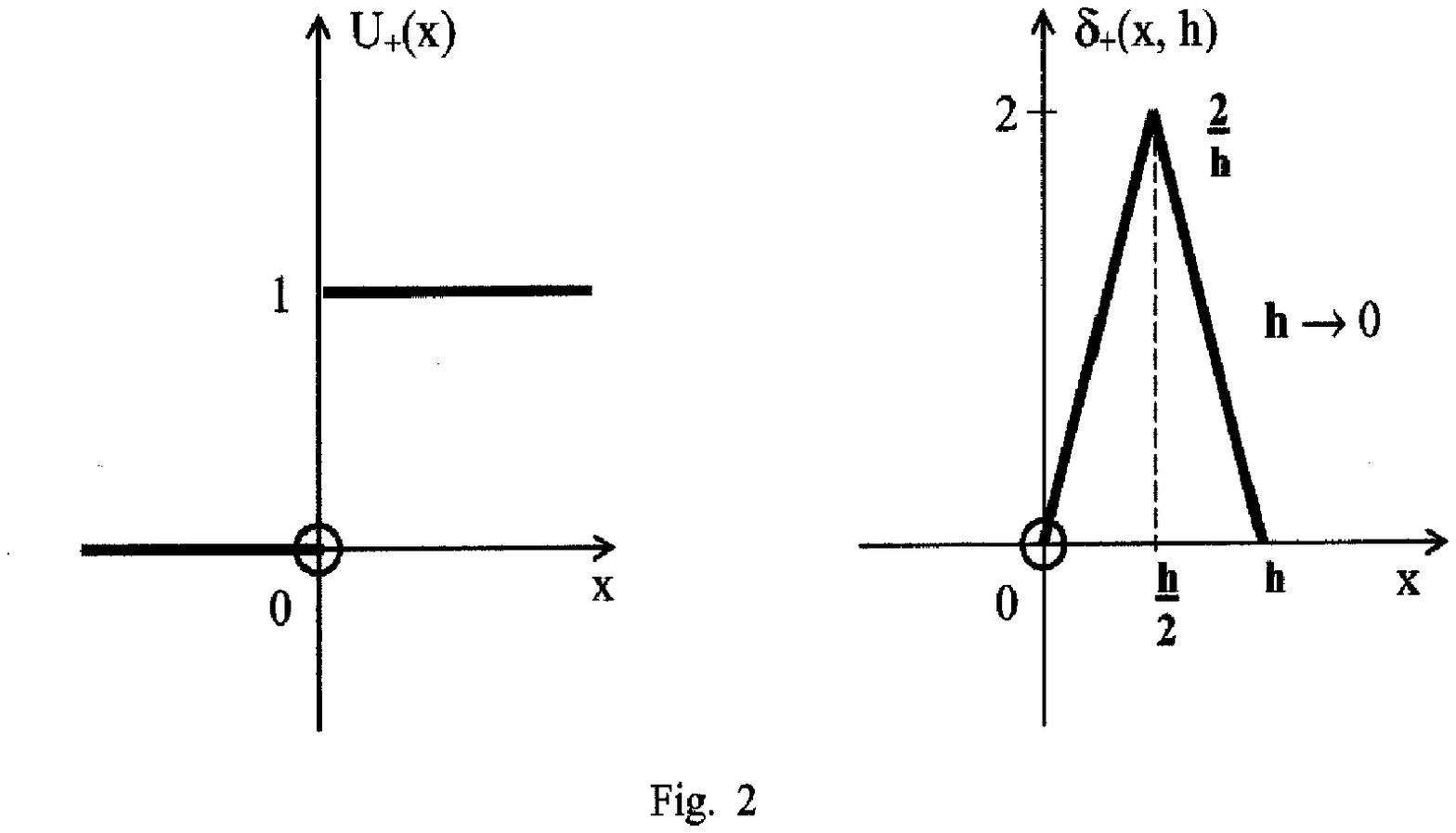,width=15cm,height=8cm}}
\medskip

\vspace{6mm}

In the papers [4], [5] and in this paper it has emerged all the
importance of the tools like the Laplace transformation and the
Stieltjes integral with certain facts from the theory of functions
of a complex variable in order to characterizing the analytic
properties of the completely monotonic functions.

\vspace{1mm}

In particular we recall that the Stieltjes integral is much used in
Mechanics and Probability,  since it unifies the treatment of the
continuous and discrete (and mixed) distributions of mass or
probability.

\vspace{2mm}

If $\alpha \left( x \right)$ is piecewise differentiable, then
$d{\kern 1pt} \alpha \left( x \right) = \alpha {\kern 1pt} '\left( x
\right){\kern 1pt} dx$, and the Stieltjes integral is simply in the
following form (reduction of a Stieltjes integral to a Riemann
integral):

\[
\int\limits_a^b {f\left( x \right)\,} \,\alpha '\left( x \right)dx
\]

For instance a real and similar case to ours  is the following: if
$\alpha \left( x \right)$ is a Heaviside step function, with point
masses $m_i $ at $x = x_i $, then

\[
d{\kern 1pt} \alpha \left( {x_i } \right) = \mathop {\lim
}\limits_{\varepsilon \downarrow 0} \left[ {\alpha {\kern 1pt}
\left( {x_i + \varepsilon } \right) - \alpha {\kern 1pt} \left( {x_i
- \varepsilon } \right)} \right]\; = m_i , \quad \int\limits_a^b
{f\left( x \right)\,} d\,\alpha \left( x \right) = \sum\limits_i
{m_i f\left( {x_i } \right)}
\]

In this case the integration by parts is usual:

\[
\int\limits_a^b {\,f\left( x \right)\,dg\left( x \right)} =
\;f\left( x \right) \cdot g\left( x \right)\left| {_a^b }
\right.\,\; - \,\,\int\limits_a^b {\,g\left( x \right)\,df\left( x
\right)}
\]

\vspace{2mm}

Suppose ${\kern 1pt} \alpha \left( 0 \right) = 0$ , ${\kern 1pt}
\alpha \left( x \right) = o\left( {e^{cx}} \right)$ as $x \to \infty
{\kern 1pt} $, and that $\Re s \ge c$, then

\vspace{2mm}

\begin{equation}
\label{eq36} \int\limits_0^\infty {e^{ - sx}} d\,\alpha \left( x
\right) = s \cdot \int\limits_0^\infty {\alpha \left( x \right)\,e^{
- sx}} dx
\end{equation}

\vspace{2mm}

The integral on the left side represents a Laplace-Stieltjes
transform, while the integral on the right side is an ordinary
Laplace transform.

\vspace{3mm}

More precisely (see [4], compare the equation 14) the second member
of (\ref{eq36}) is called s-multiplied Laplace transform.

\vspace{3mm}

However, we also recall the following result, due to S. Bernstein
([2], pp. 439-440), who was the starting point of many researches in
the Probability theory:

\vspace{3mm}

\textbf{Theorem}:\textit{ A function }$\psi \left( s \right)$
\textit{on }$\left( {0,\infty } \right)$ \textit{ is the Laplace
transform of a probability distribution }$F\left( x \right)$:

\[
\psi \left( s \right)\; = \int\limits_0^\infty {\,e^{ - sx}dF\left(
x \right)}
\]

\noindent \textit{ if and only if it is completely monotone in}
$\left( {0,\infty } \right)$ with $\psi \left( {0 + } \right) = 1$.

\vspace{3mm}

In particular, an immediate example for this beautiful theorem is
the following function:

\vspace{3mm}

\[
\psi \left( s \right) = \;\frac{\mbox{<^s>}\;(s)}{2}
\]

\vspace{3mm}

Notice that the Dirac delta function may be interpreted as a
probability density function and that the cumulative distribution
function is the Heaviside step function.

\vspace{3mm}

It is known that if $X$ is a random variable, the corresponding
probability distribution assigns to the interval [$a$, $b$] the
probability Pr[$a   \le  X   \le  b$];  for example the probability
that the variable $X$ will take a value in the interval [$a$, $b$].

\vspace{3mm}

Now, the probability distribution of the variable $X$ can be
uniquely described by its cumulative distribution function $F(x)$,
which is defined by:

\[
F\left( x \right) = \Pr \left[ {X \le x} \right]
\]

\noindent for any $x \in \mathbf{R}$ and where the right-hand side
represent the probability that the variable $X$ takes on a value
less than or equal to $x$.

\vspace{3mm}

We observe that, in our case, the Heaviside step function is the
cumulative distribution function of a random variable which is
almost surely 0.

\vspace{3mm}

\newpage

\begin{center}
\textbf{APPENDIX}
\end{center}

\vspace{4mm}

The choice of the Arabian letter <^s> $ $ is simply a consequence of
the following four reasons:
\begin{itemize}
\item  the exhaustion of Latin, Greek, Gothic, Jewish
 and other letters
to name a new special function;

\item the reference to the strokes of continuous curves (small
arcs) that characterize the same function ( [4], see Fig. 3);

\item the three dots over the letter that one by one remember
the three second-order Eulerian numbers: \textbf{1, 120, 494},
inside the Eulerian triangle ( [4], see Fig. 6);

\item the meaning of such a letter, that in Al Karaji's algebra
(about 1000 A.D.)\footnote{ \textit{Al-Karaji Abu Bakr Muhammad ibn
al-Hasan (fl. c. 400/1009) : texts and studies} / Collected and
repr. by Fuat Sezgin. - Frankfurt/M.: Institut f\"{u}r Geschichte
der Arabisch-Islamischen Wissenschaften, 1998. (Publications of the
Institute for the History of Arabic-Islamic Science : Islamic
mathematics and astronomy; n. 38). } was declared as the unknown
``par excellence'', that is absolutely comparable to our ``x'', but
also with the characteristic to form the only conjunction-``ring''
between the world of the \textit{unknown} (algebra) and the world of
the \textit{known} (arithmetics).
\end{itemize}

\vspace{3mm}

Actually it seems that the origin of the symbol  <^s> $ $, that is
pronounced \textit{shin }, is in the Ancient Egypt.

\vspace{4mm}

In fact the hieroglyph
{\includegraphics[width=0.40in,height=0.25in]{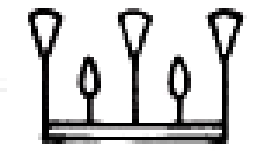}}, that
is the Egyptian syllable \textit{sha}, is similar and it represents
papyrus plants along the Nile.

 \vspace{4mm}

{\bf Acknowledgements.} This paper and the papers [4], [5] are
dedicated to the memory of my father Alessandro Ossicini.

\end{document}